\documentclass{siamltex}
\usepackage{graphicx}
\usepackage{bm}
\usepackage{amsmath}
\usepackage{amssymb}
\usepackage{epsfig}
\usepackage{subfigure}
\usepackage{latexsym}
\usepackage{amsbsy}
\usepackage{multirow}
\usepackage{lineno}
\usepackage{xparse}
\usepackage[usenames, dvipsnames]{color}
\usepackage[breaklinks,colorlinks=true,allcolors=blue]{hyperref}

\NewDocumentCommand{\dgal}{sO{}m}{%
  \IfBooleanTF{#1}
    {\dgalext{#3}}
    {\dgalx[#2]{#3}}%
}
\NewDocumentCommand{\dgalext}{m}{%
  \sbox0{%
    \mathsurround=0pt 
    $\left\{\vphantom{#1}\right.\kern-\nulldelimiterspace$%
  }%
  \sbox2{\{}%
  \ifdim\ht0=\ht2
    \{\kern-.45\wd2 \{#1\}\kern-.45\wd2 \}%
  \else
    \left\{\kern-.5\wd0\left\{#1\right\}\kern-.5\wd0\right\}%
  \fi
}

\NewDocumentCommand{\dgalx}{om}{%
  \sbox0{\mathsurround=0pt$#1\{$}%
  \sbox2{\{}%
  \ifdim\ht0=\ht2
    \{\kern-.45\wd2 \{#2\}\kern-.45\wd2 \}%
  \else
    \mathopen{#1\{\kern-.5\wd0 #1\{}
    #2
    \mathclose{#1\}\kern-.5\wd0 #1\}}
  \fi
}

\newtheorem{remark}{Remark}[section]

\newtheorem{assumption}{Assumption}[section]

\newcommand{\sigmab}{\mbox{\boldmath$\sigma$}}

\font\msbm=msbm10

\newcommand{\R}{\hbox{{\msbm \char "52}}}

\def\rbl#1{{\textcolor{black}{#1}}}

\def\rrb#1{{\textcolor{black}{#1}}}

\definecolor{otherblue}{rgb}{0,0.3,0.6}

\def\cpbl#1{{\textcolor{black}{#1}}}

\begin{document}

\bibliographystyle{plain}

\title{{Parameter-robust Stochastic Galerkin mixed approximation for linear poroelasticity with uncertain inputs}\thanks{{This work was supported 
by  EPSRC grant EP/P013317. }}}

\author{Arbaz~Khan \thanks{Department of Mathematics, Indian Institute of Technology Roorkee (IITR), Roorkee, India (\tt{arbaz@ma.iitr.ac.in})}
\and 
Catherine E. Powell\thanks{Department of Mathematics, University of Manchester, UK (\tt{c.powell@manchester.ac.uk})}
}

\pagestyle{myheadings}
\markboth{A.~Khan, C.E.~Powell}{Stochastic Galerkin Mixed FEM for Biot's Consolidation model}

\maketitle

\begin{abstract} Linear poroelasticity models have a number of important applications in biology and geophysics. In particular, Biot's consolidation model is a well-known model that describes the coupled interaction between the linear response of a porous elastic medium and a diffusive fluid flow within it, assuming small deformations. Although deterministic linear poroelasticity models and finite element methods for solving them numerically have been well studied, there is little work to date on robust algorithms for solving poroelasticity models with uncertain inputs and for performing uncertainty quantification (UQ).   The Biot model has a number of important physical parameters and inputs whose precise values are often uncertain in real world scenarios. In this work, we introduce and analyse the well-posedness of a new five-field model with uncertain and spatially varying Young's modulus and hydraulic conductivity field. By working with a properly weighted norm, we establish that the weak solution is stable with respect to variations in key physical parameters, including the Poisson ratio. We then introduce a novel locking-free stochastic Galerkin mixed finite element method that is robust in the incompressible limit. Armed with the `right' norm, we construct a parameter-robust preconditioner for the associated discrete systems. Our new method facilitates forward UQ, allowing efficient calculation of statistical quantities of interest and is provably robust with respect to variations in the Poisson ratio, the Biot--Willis constant and the storage coefficient, as well as the discretization parameters.


\end{abstract} 



\section{Introduction}\label{sec11}  Over the last three decades, there has been a substantial amount of work on developing \emph{deterministic} mathematical models of poroelasticity and finite element methods for their numerical solution. Such models have a wide range of applications in science and engineering. In particular, linear poroelasticity models have important applications in biology and medicine, such as in the study of diseases affecting the cervical spinal cord (see \cite{stoverud2016poro} and references therein) as well as applications in geophysics \cite{coussy2004poromechanics,wang2000theory}.  

The Biot consolidation model is a popular model that describes the coupled response of a linear elastic porous medium and a diffusive fluid flow within it, when subject to small deformations. The standard model is time-dependent and consists of a momentum conservation equation derived under the quasi-static assumption and a fluid mass conservation law. 
In this work, we start by considering the static two-field Biot model from \cite{lee2017parameter} which is derived by applying an implicit time discretization scheme to the standard time-dependent model.  This can be thought of as a model to be solved at a single time step. In this setting, given a body force $\bm{f}$ and a volumetric source/sink term $g$, one aims to find the displacement $\bm{u}$ of the saturated poroelastic medium (the `material') and the associated fluid pressure $p_{F}$ satisfying
\begin{subequations} \label{os2anew}
\begin{align}
-\nabla\cdot\sigmab&= \bm{f}\quad\mbox{in }D,\\
-s_0 p_F-\alpha\nabla\cdot \bm{u}+ \tau \, \nabla \cdot (\tilde{\kappa} \nabla p_F)&=g\quad \mbox{in } D,
\end{align}
\end{subequations}
with (for simplicity) homogeneous boundary conditions
\begin{subequations} \label{os2abonnew}
\begin{align}
\sigmab\bm{n}&={\bm{0}},\quad p_F=  0 \quad \qquad \qquad   \mbox{on }\partial D_{p}.\\
 \bm{u}&= \bm{0}, \quad  (\tilde{\kappa}\nabla p_F)\cdot\bm{n} =0 \, \, \quad\mbox{on }\partial D_{\bm{u}}.
\end{align}
\end{subequations}
Here, $\tau$ (with $0<\tau < < 1$) denotes the chosen time-step and the stress and strain tensors are defined as 
$$\bm{\sigma}:=2\mu\bm{\epsilon}(\bm{u})+\lambda \nabla\cdot\bm{u}\mathbf{I}-\alpha p_F\mathbf{I}, \qquad \bm{\epsilon}(\bm{u}) :=(\nabla \bm{u} +(\nabla\bm{u})^\top )/2,$$
respectively, where $\mathbf{I}$ is the $d \times d$ identity matrix (with $d=2,3$). For the analysis that follows, we assume that the spatial domain $D$ is a bounded Lipschitz polygon in $\R^2$ (polyhedral in $\R^3$) and the boundary $\partial D=\partial {D}_{\bm{u}}\cup\partial {D}_p$ is partitioned into two parts with $\partial {D}_{\bm{u}}\cap \partial D_{p}=\emptyset$ and $\partial D_{\bm{u}}, \partial D_p\neq \emptyset$.  
 
The boundary-value problem \eqref{os2anew}--\eqref{os2abonnew} has a number of important physical parameters. The Biot--Willis coefficient $\alpha \in (0,1]$, $\tilde{\kappa} >0$ is the (spatially varying) hydraulic conductivity, which depends on the permeability of the medium and viscosity of the fluid and $\lambda >0$, $\mu>0$ are the usual Lam\'{e} coefficients which can be written in terms of  the Young modulus $E>0$ and the Poisson ratio $\nu \in (0, 1/2)$ as follows
 $$\mu = \frac{E}{2(1+ \nu)}, \qquad \lambda = \frac{E \nu}{(1+ \nu)(1-2\nu)}. $$
Note that when the material becomes nearly incompressible, $\nu \to 1/2$ and $\lambda \to \infty$.  The parameter $s_0$ is the so-called storage coefficient and for incompressible fluids we typically have
 \begin{align}\label{s0_def}
 s_0=\frac{\alpha-\phi}{2\mu d^{-1}+\lambda},
 \end{align}
 where $\phi$ denotes porosity. Using the definition of $\alpha$, we assume $0 \le \phi < \alpha \le 1$ and so $s_{0} > 0$. For other characterisations of the storage coefficient see  \cite{detournay1993fundamentals,merxhani2016introduction,lee2017parameter,botti2019numerical}.  We also refer interested readers to \cite{wang2000theory,coussy2004poromechanics} for more background on poromechanics. 
 
 
 

A broad discussion on the use of finite element methods (FEMs) for the numerical simulation of fluid flow and deformation processes in porous media can be found in \cite{lewis1998finite}. Taylor-Hood finite element approximation for Biot consolidation models is discussed in particular in \cite{murad1992improved,murad1994stability, murad1996asymptotic} and a stabilized lowest-order FEM for a three-field model is discussed in \cite{berger2015stabilized}. It is well accepted that solving poroelasticity models such as \eqref{os2anew}--\eqref{os2abonnew} numerically is challenging. The main issue  (e.g., see \cite{phillips2009overcoming,oyarzua2016locking,haga2012causes}) is that the accuracy of standard FEMs deteriorates due to spurious pressure modes and volumetric locking when $\nu \to 1/2$.  To avoid locking, various \emph{mixed} formulations have been proposed which are derived by introducing appropriate auxiliary variables. In \cite{oyarzua2016locking}, Oyarz\'{u}a et al.~discussed locking-free FEMs and presented a priori error analysis for a three-field mixed Biot consolidation model. In \cite{lee2017parameter}, Lee et al.~presented an alternative three-field model and discussed parameter-robust mixed FEM approximation.  The authors of \cite{lee2019mixed} also recently presented a mixed method for a nearly incompressible multiple-network poroelasticity model. In addition to developing methods that avoid locking, a second important challenge is to design bespoke solvers for the associated discrete linear systems that are robust with respect to both variations in the discretization parameters \emph{and} the physical parameters. See \cite {lee2017parameter} and \cite{relax_prec} for recent work on preconditioning in the context of poroelasticity.

Adopting a similar approach as in \cite{lee2017parameter} and \cite{oyarzua2016locking},  we introduce the `total pressure' $p_T:=-\lambda \nabla \cdot\bm{u}+\alpha p_F$ and consider the following three-field `mixed' model: 
\begin{subequations} \label{os2a}
\begin{align}
-\nabla\cdot\sigmab&= \bm{f}\quad\mbox{in }D,\\
-\nabla\cdot \bm{u}-\lambda^{-1}(p_T-\alpha p_F)&=0\quad \mbox{in }D,\\
\lambda^{-1}(\alpha p_T-\alpha^2 p_F)-s_0 p_F+\nabla \cdot (\kappa \nabla p_F)&=g\quad \mbox{in } D.
\end{align}
\end{subequations}
Here, the stress tensor is defined as $\bm{\sigma}:=2\mu\bm{\epsilon}(\bm{u})-p_T\mathbf{I}$ and the rescaled hydraulic conductivity is $\kappa:=\tau \tilde{\kappa}. $ Note that when $\nu \to 1/2$ we have $\lambda^{-1} \to 0$ and $s_{0} \to 0$. The problem then decouples for $(\bm{u}, p_{T})$ and $p_{F}$ and the weak solution remains well behaved.  If we set $s_0=\alpha^2/\lambda$ in \eqref{os2a} then we recover the model in \cite{lee2017parameter}. The model in \cite{oyarzua2016locking} differs in that it is derived from a two-field model with a slightly different stress tensor.  
 
 



In real-world applications, the values of $E, \kappa, \nu, \alpha$ and $s_{0}$ may have very different orders of magnitude. As reported in \cite{lee2017parameter}, in biomedical applications involving fluid flow in the soft tissue of the central nervous system (e.g., see \cite{stoverud2016poro}), $E$ typically takes values in the range $1-60$ kPA, while $\nu$ takes values from $0.3$ to almost $0.5$ (the incompressible limit)  and the permeability lies in the range $10^{-14}$--$10^{-16}\, m^2$. In geophysical applications (e.g., see \cite{coussy2004poromechanics,wang2000theory}), $E$ can be in the order of GPA,  while  $\nu$ varies from $0.1$ to almost $0.3$ and the permeability is in the range $10^{-9}$--$10^{-21}\,m^2$. While realistic ranges of values for the inputs may be available, their precise values are often \emph{uncertain}. Even if measurements are available, these are often subject to errors. Moreover, materials may have small imperfections or variations which are impossible to characterise, so quantities such as  $E$ and $\kappa$ may be spatially varying in an uncertain way.  


Although there is a large body of work on deterministic poroelastic models, there has been little work to date on robust approximation for \emph{stochastic} Biot models. That is, formulations in which one or more inputs is modelled as a function of random variables (or parameters). In \cite{chang1985uncertainty}, uncertainty in one-dimensional consolidation of soils was assessed using the method of moments and Monte Carlo (MC) simulation.  Uncertainty in consolidation of soils was also considered in \cite{nishimura2002consolidation} and \cite{darrag1993consolidation}.  In \cite{frias2004stochastic}, Frias et al.~discussed stochastic modelling of highly heterogeneous poroelastic media with long-range correlations using finite element approximation and MC sampling. In \cite{delgado2015stochastic}, Delgado et al.~outlined a stochastic Galerkin finite element method (SGFEM) for a two-field poroelastic model and demonstrated its use on a model with spatially uniform random inputs. More recently, Botti et al.~\cite{botti2019numerical} discussed the numerical solution of a two-field Biot model with random coefficients using a non-intrusive polynomial chaos method. To the best of our knowledge, there are no previous works addressing parameter-robust \emph{mixed} formulations of stochastic Biot consolidation models that are suitable in the nearly incompressible case. To remedy this, we combine and extend ideas from  \cite{lee2017parameter} and \cite{oyarzua2016locking} (for deterministic Biot models) and \cite{KPSpre} (for a stochastic linear elasticity model) and formulate, analyse and then solve a new mixed formulation of Biot's consolidation model with uncertain inputs. We tackle the case where $E$ and $\kappa$ are random fields with prescribed distributions and employ mixed SGFEM approximation. Our new method facilitates efficient forward uncertainty quantification (UQ), allowing calculation of statistical quantities of interest and is provably robust with respect to variations in $\nu, \alpha$ and (after rescaling) to $s_{0}$. 

\section{The New Model}
To define the new model, we introduce vectors of parameters $\bm{y}=(y_1,\ldots,y_{M_1})$ and $\bm{z}=(z_1,\dots,z_{M_2})$ which are assumed to be images of mean-zero, bounded and independent real-valued random variables.  We consider models where $E$ and $\kappa$ are expressed as functions of the form
\begin{align}
\label{E_def} E(\bm{x},\bm{y})&= e_0(\bm{x})+\sum_{k=1}^{M_1}e_k(\bm{x}) y_k, \quad \bm{x}\in D, \, \bm{y}\in\Gamma_y:=\Gamma_1\times\cdots\times\Gamma_{M_1},\\
\label{k_def} \kappa(\bm{x},\bm{z})&= \kappa_0(\bm{x})+\sum_{k=1}^{M_2}\kappa_k(\bm{x}) z_k, \quad\bm{x}\in D,  \, \bm{z}\in\Gamma_z:=\Gamma_1\times\cdots\times\Gamma_{M_2}.
\end{align}
For simplicity, we further assume that $y_k,z_k\in\Gamma_k=[-1,1]$ for each $k$ (but any bounded interval is permitted for the analysis that follows). Now, if we define the parameter domain $\Gamma:=\Gamma_y\times\Gamma_z$, and let $\bm{Y}:=(\bm{y},\bm{z})$, the parametric analogue of \eqref{os2a} is:  find $\bm{u}: D\times\Gamma \to \rrb{\mathbb{R}^{d}}$ and ${{p}}_T, {p}_F:D\times\Gamma \to \mathbb{R}$ such that 
\begin{subequations} \label{sgos2a}
\begin{align}
-\nabla\cdot\sigmab(\bm{x},\bm{Y})&= \bm{f}(\bm{x})\quad\mbox{in }D\times \Gamma, \\
-\nabla\cdot \bm{u}(\bm{x},\bm{Y})-\frac{(p_T(\bm{x},\bm{Y})-\alpha p_F(\bm{x},\bm{Y}))}{\lambda(\bm{x},\bm{y})}&=0 \quad \, \, \, \, \quad \mbox{in }D\times \Gamma, \\
\frac{\alpha p_T(\bm{x},\bm{Y})}{\lambda(\bm{x},\bm{y})}- \left(\frac{\alpha^2}{\lambda(\bm{x},\bm{y})} +s_0(\bm{x},\bm{y})\right) p_F(\bm{x},\bm{Y})+\nabla \cdot (\kappa(\bm{x},\bm{z}) \nabla p_F(\bm{x},\bm{Y}))&=g(\bm{x}) \, \, \quad \mbox{in } D\times\Gamma.
\end{align}
\end{subequations}
Note that \eqref{E_def} and \eqref{k_def} have the same structure as truncated Karhunen--Lo\`eve expansions. In that setting $e_{0}$ and $k_{0}$ represent the means of the associated random fields. To distinguish them from the \emph{physical} parameters, we will refer to $\bm{Y}=(\bm{y}, \bm{z})$ as the \emph{stochastic} parameters.  In \eqref{sgos2a}, the solution fields $\bm{u}, p$ and $\tilde{p}$ are all functions of $\bm{Y}$, as are the stress  and strain tensors $\sigmab : D\times\Gamma\rightarrow \R^{d\times d}$ and $\bm{\epsilon} : D\times\Gamma\rightarrow \R^{d\times d}$. Since they depend on $E$, the Lam\'{e} coefficients are now functions of the stochastic parameters $\bm{y}$. That is,
\begin{align*}
\mu(\bm{x},\bm{y})=\frac{E(\bm{x},\bm{y})}{2(1+\nu)}, \quad \lambda(\bm{x},\bm{y})=\frac{E(\bm{x},\bm{y})\nu}{(1+\nu)(1-2\nu)}, \qquad \bm{x} \in D, \, \bm{y} \in \Gamma_{y}.
\end{align*}

Stochastic Galerkin (SG) approximation is one of the most well-known approaches for performing forward UQ  in parametric PDEs.  Unlike many of its non-intrusive competitors, it offers a natural framework for error analysis and a posteriori error estimation \cite{khan2018error3field}. From a practical perspective, \emph{standard} SGFEMs (which employ finite element approximation for the spatial discretisation) can be applied straightforwardly if  (i) the model has a modest number $M$ of stochastic parameters, (ii)  the PDE inputs are expressed as \emph{linear} functions of the stochastic parameters as in \eqref{E_def}--\eqref{k_def} and (iii) efficient linear algebra tools are available. SGFEMs also have excellent convergence properties for models whose solutions are smooth functions \cite{Cohen} of the stochastic parameters (such as scalar elliptic PDEs with affine parameter dependence). However, for more challenging problems, more advanced adaptive and/or multilevel SGFEMs may be required. Such methods are now available for several classes of parametric PDEs (see \cite{Adam_ML}, \cite{EIGEL}, \cite{khan2018error3field}). See also \cite{schwab2011sparse} for a comprehensive discussion of SGFEMs for high-dimensional parametric PDEs.



Recall that $E$ and $\kappa$ have the form \eqref{E_def} and \eqref{k_def}, respectively. Unfortunately, $E^{-1}$ also appears in three places in \eqref{sgos2a} (due to the $1/\lambda$ term which was introduced via $p_{T}$ in setting up the three-field mixed model to accommodate the nearly incompressible case).  Clearly, $E^{-1}$ is not linear in the stochastic parameters.  While it is not impossible to apply a stochastic Galerkin method directly to \eqref{sgos2a}, the cost of assembling and solving the associated linear system can be problematic for stochastically nonlinear problems.  Following ideas in  \cite{KPSpre} and \cite{khan2018error3field}, we reformulate the model so that the resulting discrete problem, although larger, has a more favourable structure and a sparser coefficient matrix.  To this end, we introduce two auxiliary variables ${p}_1:=(p_T-\alpha p_F)/E,$ $p_2:= p_F/E,$ and define the rescaled Lam\'e and storage coefficients
\begin{align}\label{new_constants}
\tilde{\mu} := \frac{2\mu}{E} = {1\over 1+ \nu}, \qquad \tilde{\lambda} := \frac{\lambda}{E} =\frac{\nu}{(1+\nu)(1-2\nu)}, \qquad \tilde{s}_0:= Es_{0}.
\end{align}
Note that if $s_{0}$ is defined as in \eqref{s0_def} then $\tilde{s}_{0}$ is independent of $E$.

Now, substituting $p_{1}$ and $p_{2}$ in \eqref{sgos2a} and rearranging yields the following five-field formulation:  find $\bm{u}: D\times\Gamma \to \rrb{\mathbb{R}^{d}}$ and ${{p}}_T, {p}_F,p_1,p_2:D\times\Gamma \to \mathbb{R}$ such that, 
\begin{subequations} \label{sgos2a11}
\begin{align}
-\nabla\cdot\sigmab(\bm{x},\bm{Y})&= \bm{f}(\bm{x})\quad\mbox{in }D\times \Gamma, \\
-\nabla\cdot \bm{u}(\bm{x},\bm{Y})- \tilde{\lambda}^{-1} p_{1}(\bm{x},\bm{Y})
&=0 \quad \, \, \, \quad \mbox{in }D\times \Gamma, \\
\tilde{\lambda}^{-1}\alpha \, p_1(\bm{x},\bm{Y})- \tilde{s}_0 \, p_2(\bm{x},\bm{Y})+\nabla \cdot (\kappa(\bm{x},\bm{z}) \nabla p_F(\bm{x},\bm{Y}))&=g(\bm{x})\quad \mbox{in } D\times\Gamma,  \\
-\tilde{\lambda}^{-1}(p_T(\bm{x},\bm{Y})-\alpha p_F(\bm{x},\bm{Y}))+ \tilde{\lambda}^{-1} E(\bm{x},\bm{y})p_1(\bm{x},\bm{Y})&=0 \quad \, \, \quad\mbox{in } D\times \Gamma, \\
-\tilde{s}_0 p_F(\bm{x},\bm{Y})+ \tilde{s}_0E(\bm{x},\bm{y})p_2(\bm{x},\bm{Y})&=0 \quad \, \, \quad\mbox{in } D\times \Gamma.
\end{align}
\end{subequations}
Notice that $E$ appears in the first, fourth and fifth equations and $\kappa$ appears in the third equation, but $E^{-1}$ does not appear at all. We also supplement \eqref{sgos2a11} with the following boundary conditions
\begin{subequations} \label{sgos2a11B}
\begin{align}
\sigmab\bm{n}&={\bm{0}}, \quad p_F=  0 \quad \qquad \qquad   \mbox{on }  \partial D_p \times \Gamma\\
\bm{u}&= \bm{0}, \quad (\kappa \nabla p_F)\cdot\bm{n} =0 \quad  \, \, \mbox{on }  \partial D_{\bm{u}} \times \Gamma.
\end{align}
\end{subequations}


\subsection{Outline} In Section \ref{EPWRD} we develop the weak formulation of  \eqref{sgos2a11}--\eqref{sgos2a11B} and prove that it is well-posed.  In Section \ref{disformul},  we derive a finite-dimensional weak problem using a stochastic Galerkin mixed FEM (SG-MFEM) and discuss the block structure of the associated discrete linear system and coefficient matrix. In Section \ref{precond-sec}, we apply the theory of operator preconditioning outlined in \cite{lee2017parameter} and build on previous work \cite{KPSpre} for linear elasticity problems, to develop a new block-diagonal preconditioner that is provably robust with respect to $\nu$, $\alpha$ and $\tilde{s}_{0}$ (the Poisson ratio, the Biot--Willis constant and the rescaled storage coefficient). Finally, in Section \ref{Numer}, we present numerical results to demonstrate the robustness and efficiency of the proposed solver and then use the new SG-MFEM scheme to perform forward UQ and estimate statistical information about the displacement $\bm{u}$ and fluid pressure $p_{F}$ for a benchmark `footing' problem.



\section{Weak Mixed Formulation}\label{EPWRD}
 We begin by introducing some function spaces  for the forthcoming analysis and state some assumptions on $E$ and $\kappa$ defined in \eqref{E_def}--\eqref{k_def}. Recall that $y_{k}, z_{k} \in [-1,1]$, $\bm{y} \in \Gamma_{y}:=[-1,1]^{M_{1}}$ and $\bm{z} \in \Gamma_{z}:=[-1,1]^{M_{2}}$. In addition,  $\bm{Y}=(\bm{y}, \bm{z}) \in \Gamma=[-1,1]^{M}$ where $M:=M_{1}+M_{2}$. 
 
 \smallskip 
 
 \begin{assumption}\label{Assump2}
$E\in L^{\infty}(D\times \Gamma_y)$ and $\kappa\in L^{\infty}(D\times \Gamma_z)$ and there exist positive constants $E_{\min}$, $E_{\max}$, $\kappa_{\min}$ and $\kappa_{\max}$ such that
\begin{align}\label{bounde11}
&0<E_{\min}\le E(\bm{x},\bm{y}) \le E_{\max} <\infty \quad \mbox{\rm{a.e. in}}\, D\times \Gamma_y,\\
&0<\kappa_{\min}\le \kappa(\bm{x},\bm{z}) \le \kappa_{\max} <\infty \;\;\quad \mbox{\rm{a.e. in}}\, D\times \Gamma_z.
\end{align}
In addition, there exist positive constants $e_0^{\min}$, $e_0^{\max}$, $\kappa_0^{\min}$ and $\kappa_0^{\max}$ such that
\begin{align}
\label{bounde11x}
\cpbl{0<e_0^{\min}\le e_0(\bm{x})\le e_0^{\max} < \infty  \quad \mbox{\rm{a.e. in}}\, D \quad \mbox{and}  \quad  \sum_{k=1}^{M_1}|| e_{k} ||_{L^{\infty}(D)}< e_0^{\min}},\\
\cpbl{0<\kappa_0^{\min}\le \kappa_0(\bm{x})\le \kappa_0^{\max} < \infty  \quad \mbox{\rm{a.e. in}}\, D \quad \mbox{and}  \quad  \sum_{k=1}^{M_2}|| \kappa_{k} ||_{L^{\infty}(D)}< \kappa_0^{\min}}.
\end{align}
\end{assumption}

Next, we define a product\footnote{This is natural if $Y_{k},$ $k=1,\ldots, M$ are associated with an underlying set of \emph{independent} random variables.} measure $\pi(\bm{Y}):= \Pi_{k=1}^{M} \pi_{k}(Y_{k}) = \Pi_{k=1}^{M_1}  {\pi_k(y_{k})} \, \Pi_{k=1}^{M_2}{\pi_{k+M_{1}}(z_{k})},$ 
where $\pi_k$ denotes a measure on $(\Gamma_k, \mathcal{B}(\Gamma_k))$ and $\mathcal{B}(\Gamma_k)$ is the Borel $\sigma$-algebra on $\Gamma_k=[-1,1]$. Using this measure, we can define Bochner spaces of the form,
\begin{align*}
L^2_{\pi}(\Gamma, X(D)) :=\left \{v(\bm{x},\bm{Y}) : D\times\Gamma\rightarrow \mathbb{R}; ||v||_{L^2_{\pi}(\Gamma, X(D))} <\infty\right\},
\end{align*}
where $X(D)$ is a normed vector space of real-valued functions on $D$ with norm $||\cdot||_X$ and 
\begin{align}\label{norm_def}
||\cdot||_{L^2_\pi(\Gamma,X(D))} := \left(\int_{\Gamma}||\cdot||_X^2 d\pi(\bm{Y})\right)^{1/2}.
\end{align} 
In particular, we will need the spaces $\mathcal{V}_{0,p}:= {L^2_{\pi}(\Gamma,{H}^1_{0,p}(D))}$ and  $\mathcal{W} :=  {L^2_{\pi}(\Gamma, {L}^2(D))}$ where ${H}^1_{0,p}(D)= \{ {v} \in {H}^1(D),  {v} |_{\partial D_{p}} = 0\}$  and ${H}^1(D)$ is  {the usual} Sobolev space with norm $||\cdot||_{1}$. We will also need the analogous spaces of \emph{vector-valued} functions
\begin{align*}
\bm{\mathcal{V}_{0,\bm{u}}} := {L^2_{\pi}(\Gamma,  \bm{H}^1_{0,\bm{u}}(D))}, \qquad  \bm{\mathcal{W}} :=  {L^2_{\pi}(\Gamma, \bm{L}^2(D))}
\end{align*}
where $\bm{H}^1_{0, \bm{u}}(D) = \{ \bm{v} \in\bm{H}^1(D),  \bm{v} |_{\partial D_{\bm{u}}}=\bm{0}\}$ and $\bm{H}^1(D)=\bm{H}^1(D;\mathbb{R}^{\rbl{d}})$.  Notice that $\bm{\mathcal{V}}_{0,\bm{u}}$ and $ \mathcal{V}_{0,p}$ encode the essential Dirichlet boundary conditions associated with the displacement $\bm{u}$ and fluid pressure $p_{F}$.

Assuming that the load function $\bm{f}\in \bm{L}^2(D)$ and the source term $g \in L^2(D)$, the weak form of  \eqref{sgos2a11}--\eqref{sgos2a11B} can be written as: find $(\bm{u},p_T,{p}_F,p_1,p_2)\in \bm{\mathcal{V}}_{0,\bm{u}}\times\mathcal{W}\times\mathcal{V}_{0,p}\times \mathcal{W} \times \mathcal{W}$ such that
\begin{subequations} \label{scm11a}
\begin{align}
a(\bm{u},\bm{v})+b(\bm{v},p_T)&=f(\bm{v})  \qquad \, \, \forall \bm{v}\in \bm{\mathcal{V}}_{0,\bm{u}},\\
b(\bm{u},q_T)-c_1({p}_1,q_T)&=0   \qquad \qquad  \forall q_T\in \mathcal{W},\\
c_2(p_1,{q}_F)-c_3(p_2,q_F)-d({p}_F,{q}_F)&=g(q_F)   \qquad \forall  {q}_F\in \mathcal{V}_{0,p},\\
-c_1(p_T,q_1)+c_2(p_F,q_1)+\tilde{c}_1(p_1,q_1)&=0   \qquad \qquad \forall q_1\in \mathcal{W}, \\
-c_3(p_F,q_2)+\tilde{c}_2(p_2,q_2)&=0  \qquad \qquad  \forall q_2\in\mathcal{W}.
\end{align}
\end{subequations}

In \eqref{scm11a}, the symmetric bilinear forms $a(\cdot, \cdot): \bm{\mathcal{V}}_{0,\bm{u}} \times  \bm{\mathcal{V}}_{0,\bm{u}}  \to \mathbb{R}$, $d(\cdot, \cdot): \mathcal{V}_{0,p}\times \mathcal{V}_{0,p}\to \mathbb{R}$, and $\tilde{c}_1(\cdot, \cdot), \tilde{c}_{2}(\cdot, \cdot): {\cal W} \times {\cal W} \to \mathbb{R}$ are defined by
\begin{align*} 
a(\bm{u},\bm{v})&:= \tilde{\mu}\int_{\Gamma}\int_{D}  E(\bm{x},\bm{y}) \bm{\epsilon}(\bm{u}(\bm{x},\bm{Y})):\bm{\epsilon}(\bm{v}(\bm{x},\bm{Y})) \, d\bm{x} \, d\pi(\bm{Y}),\\
d(r,s)&:=\int_{\Gamma}\int_{D} \kappa(\bm{x},\bm{z}) \nabla{r}(\bm{x},\bm{Y})\cdot\nabla{s}(\bm{x},\bm{Y}) \, d\bm{x} \, d\pi(\bm{Y}),
\end{align*}
and $\tilde{c}_1(p,q) :=\tilde{\lambda}^{-1}\tilde{c}(p,q)$,  $ \tilde{c}_2(p,q) := \tilde{s}_0 \,  \tilde{c}(p,q),$ where 
\begin{align*} 
\tilde{c}(p,q) & :=\int_{\Gamma}\int_{D} E(\bm{x},\bm{y}) {p}(\bm{x},\bm{Y}){q}(\bm{x},\bm{Y}) \, d\bm{x} \, d\pi(\bm{Y}).
\end{align*}
Recall that $\tilde{\mu}$, $\tilde{\lambda}$ and $\tilde{s}_{0}$ were defined in \eqref{new_constants}. As $\nu \to \frac{1}{2}$, $\tilde{\mu}$ remains bounded and $ \tilde{\lambda}^{-1}, \tilde{s}_{0} \to 0$.  

The bilinear form $b(\cdot, \cdot): \bm{\mathcal{V}}_{0,\bm{u}} \times \mathcal{W}\to \mathbb{R}$ is defined by
\begin{align*}
b(\bm{v},p)&:=-\int_{\Gamma}\int_{D}  {p}(\bm{x},\bm{Y}){\rm{div}}\,\bm{v}(\bm{x},\bm{Y}) \, d\bm{x} \, d\pi(\bm{Y}),
\end{align*}
and the symmetric bilinear forms $c_{1}(\cdot, \cdot), c_2(\cdot, \cdot), c_{3}(\cdot, \cdot): {\cal W} \times {\cal W} \to \mathbb{R}$ are defined by $c_1(p,q)  :=\tilde{\lambda}^{-1}c(p,q)$, $c_2(p,s) := \alpha \, \tilde{\lambda}^{-1} \, c(p,q)$, and $c_3(p,s):=\tilde{s}_0 c(p,q)$, where 
\begin{align*}
c(p,q) & :=\int_{\Gamma}\int_{D} {p}(\bm{x},\bm{Y}){q}(\bm{x},\bm{Y}) \, d\bm{x} \, d\pi(\bm{Y}).
\end{align*}
$E$ and $\kappa$ \emph{do not} appear in these expressions. Note also that although we have defined $c_{2}(\cdot, \cdot)$ and $c_{3}(\cdot, \cdot)$ as bilinear forms on ${\cal {W}} \times {\cal W}$ (for convenience), in \eqref{scm11a} they act on a pair of functions from ${\cal V}_{0,p}$ and ${\cal W}$, where ${\cal V}_{0,p} \subset {\cal W}$. Finally, we define the linear functionals $f: \bm{\mathcal{V}}_{0,\bm{u}} \to \mathbb{R}$ and $g: {\cal V}_{0,p} \to \mathbb{R}$ by, 
\begin{align*}
f(\bm{v})& : = \int_{\Gamma}\int_{D} \bm{f}(\bm{x}) \cdot \bm{v}(\bm{x},\bm{Y}) \, d\bm{x} \, d\pi(\bm{Y}), \qquad g(r) : = \int_{\Gamma}\int_{D} g(\bm{x}) {r}(\bm{x},\bm{Y}) \, d\bm{x}\, d\pi(\bm{Y}).
\end{align*}


If we now define $\bm{\mathcal{X}}:=\bm{\mathcal{V}}_{0,\bm{u}}\times\mathcal{W}\times\mathcal{V}_{0,p}\times \mathcal{W} \times \mathcal{W}$ and $B(\cdot, \cdot):\bm{\mathcal{X}} \times \bm{\mathcal{X}} \to \mathbb{R}$ by
\begin{align}\label{big_B}
B(\bm{u},p_T,{p}_F,p_1,p_2; \bm{v},q_T,{q}_F,q_1,q_2)&:=a(\bm{u},\bm{v})+b(\bm{v},p_T)+b(\bm{u},q_T)-c_1({p}_1,q_T) + c_2(p_1,{q}_F) \nonumber\\
&\quad-c_3({p}_2,{q}_F)-d(p_F,q_F)
-c_1(p_T,q_1)+c_2(p_F,q_1)\nonumber\\
&\quad+\tilde{c}_1(p_1,q_1)-c_3(p_F,q_2)+\tilde{c}_2(p_2,q_2),
\end{align}
then we can express (\ref{scm11a}) more concisely as: find $(\bm{u},p_T,{p}_F,p_1,p_2)\in \bm{\mathcal{X}}$ such that 
\begin{align}\label{scm12}
B(\bm{u},p_T,{p}_F,p_1,p_2; \bm{v},q_T,{q}_F,q_1,q_2)=f(\bm{v})+g(q_F), \quad \forall \, (\bm{v},q_T,{q}_F,q_1,q_2,q_3)\in \bm{\mathcal{X}}.
\end{align}
To establish that \eqref{scm12} is well-posed, we will need to choose an appropriate norm on $\bm{\mathcal{X}}$.  First, let $\| \cdot \|_{{\mathcal{W}}}$ and $\| \cdot \|_{\bm{\mathcal{W}}}$ denote the norms associated with ${\cal{W}}$ and $ \bm{\mathcal{W}}$ defined by \eqref{norm_def}, respectively, and note that if Assumption \ref{bounde11} holds, then $\| e_{0}^{1/2} \nabla \bm{u} \|_{\bm{\mathcal{W}}}$, $\| \kappa_{0}^{1/2} \nabla p \|_{{\mathcal{W}}}$ and  $\| e_{0}^{1/2}  p \|_{{\mathcal{W}}}$ are norms on $\bm{\mathcal{V}}_{0,\bm{u}}$, ${\mathcal{V}}_{0,p}$  and ${\mathcal{W}}$, respectively.  On $\bm{\mathcal{X}}$ we will work with the weighted and coefficient-dependent norm $||| \cdot |||$ defined by
\begin{align}\label{normden}
|||(\bm{v},q_T,q_F,q_1,q_2)|||^2 &:= \rbl{\tilde{\mu}}||e_0^{\frac{1}{2}}\nabla\bm{v}||_{\bm{\mathcal{W}}}^2+\Big(\tilde{\mu}^{-1}+{\tilde{\lambda}}^{-1}\Big)||e_0^{-\frac{1}{2}}q_T||_{\mathcal{W}}^2+ (\alpha^2 \tilde{\lambda}^{-1}+\tilde{s}_0) ||e_0^{-\frac{1}{2}}{q}_F||_{\mathcal{W}}^2\nonumber\\
&+ ||\kappa_0^{1/2}\nabla q_F||_{\mathcal{W}}^2+ \tilde{\lambda}^{-1}||e_0^{\frac{1}{2}}q_1||_{\mathcal{W}}^2+ \tilde{s}_0 ||e_0^{\frac{1}{2}}q_2||^2_{\mathcal{W}}.
\end{align}
Recall that $e_{0}$ and $\kappa_{0}$ are the leading (deterministic) coefficients in the expressions for the Young modulus $E$ and hydraulic conductivity $\kappa$, respectively.  Recall also from \eqref{new_constants} that $\tilde{\mu}$ and $\tilde{\lambda}$ depend only on the Poisson ratio $\nu$. If $s_{0}$ is defined as in \eqref{s0_def}, then $\tilde{s}_{0}$ depends on $\alpha$ and $\phi$ as well as $\nu$. 

It is straightforward to show, using Assumption \ref{bounde11} and results such as the Cauchy-Schwarz inequality, that all the bilinear forms appearing in \eqref{scm11a} are bounded in the chosen weighted norms. In particular,
\begin{align}
a(\bm{u},\bm{v}) & \le \frac{{E}_{\max}}{e_0^{\min}}(\tilde{\mu}^{1/2}||e_0^{1/2}\nabla\bm{u}||_{\bm{\mathcal{W}}} )
(\tilde{\mu}^{1/2}||e_0^{1/2}\nabla\bm{v}||_{\bm{\mathcal{W}}}) \quad\forall \, \bm{u}, \bm{v} \in \bm{\mathcal{V}}_{0,\bm{u}},  \label{a_UP_bound}\\
b(\bm{u}, p)&\le \rbl{\sqrt{d}} \, (\tilde{\mu}^{1/2}||e_0^{1/2}\nabla\bm{u}||_{\bm{\mathcal{W}}})(\tilde{\mu}^{-1/2}||e_0^{-1/2}p||_{\mathcal{W}})
\qquad  \forall \bm{u} \in \bm{\mathcal{V}}_{0,\bm{u}}, \, \forall \, p \in \mathcal{W}, \\
\tilde{c}_1(p,q)&\le \frac{E_{\max}}{e_{0}^{\min}}(\tilde{\lambda}^{-1/2} ||e_0^{1/2}p||_{\mathcal{W}}) (\tilde{\lambda}^{-1/2}||e_0^{1/2}q||_{\mathcal{W}}) \quad \,  \, \quad\forall \, p,q\in\mathcal{W}, \\
\tilde{c}_2(p,q) & \le \frac{E_{\max}}{e_{0}^{\min}}(\tilde{s}^{1/2}_0 ||e_0^{1/2}p||_{\mathcal{W}}) (\tilde{s}^{1/2}_0||e_0^{1/2}q||_{\mathcal{W}}) \, \quad \qquad \, \, \, \, \forall  \, p,q\in\mathcal{W},\\
d(p,q)&\le \frac{\kappa_{\max}}{\kappa_{0}^{\min}}||\kappa_0^{1/2}\nabla p||_{\mathcal{W}} ||\kappa_0^{1/2}\nabla q||_{\mathcal{W}} \quad \quad \qquad \qquad \quad \, \, \forall \, p,q\in\mathcal{V}.
\label{dbound}
\end{align}
In addition, we have the following coercivity results,
\begin{align}
a(\bm{u},\bm{u})& \ge \frac{{E}_{\min}}{e_{0}^{\max}} 
C_{K} \, \tilde{\mu}||e_0^{1/2}\nabla\bm{u}||_{\bm{\mathcal{W}}}^2 \quad \, \,  \forall \, \bm{u}\in \bm{\mathcal{V}}_{0,\bm{u}}, 
\label{abound}\\
\tilde{c}_1(p,p)&\ge\frac{E_{\min}}{e_0^{\max}}  \, \tilde{\lambda}^{-1} ||e_0^{1/2}p||^2_{\mathcal{W}} \quad \quad \, \quad\forall  \, p\in \mathcal{W},
\label{c1bound} \\
\tilde{c}_2(p,p) &  \ge\frac{E_{\min}}{e_0^{\max}} \, \tilde{s}_0 ||e_0^{1/2} p||^2_{\mathcal{W}}  \qquad \quad  \, \quad\forall \,  p\in \mathcal{W},
\label{c2bound} \\
d(p,p)&\ge \frac{\kappa_{\min}}{\kappa_{0}^{\max}}\, ||\kappa_0^{1/2} \nabla p||^2_{\mathcal{W}} \quad \qquad \quad \forall  \, p\in \mathcal{W},
\label{dubound}
\end{align}
where $\rbl{0<C_K \leq 1}$ is the usual Korn constant (see \cite{KPSpre} and references therein). Finally, following Lemma 2.2 from \cite{KPSpre}, it can also be shown that there exists a constant $C_{D}>0$ (the inf-sup constant) such that
\begin{align}\label{first_inf_sup}
\sup_{0\neq\bm{v}\in \bm{\mathcal{V}}_{0,\bm{u}} }
\frac{b(\bm{v},q)}{ ||\nabla\bm{v}||_{\bm{\mathcal{W}}} }&\ge C_{D} 
 ||q||_{\mathcal{W}} \quad\forall  q\in \mathcal{W}.
 \end{align}
Combining all these bounds, the well-posedness of \eqref{scm12} is established using the next result.




\smallskip

\begin{lemma}\label{main_lemma} Suppose $\alpha^2\tilde{\lambda}^{-1}\in(0,\tilde{C}_1\tilde{s}_0]$ with $\tilde{C}_1<\frac{3}{2}$. For  any $(\bm{u},p_T,{p}_F,p_1,p_2) \in \bm{\mathcal{X}}$, there exists  $(\bm{v},q_T,{q}_F, {q}_1, {q}_2) \in \bm{\mathcal{X}}$ with 
$ ||| (\bm{v},q_T,{q}_F, {q}_1, {q}_2) |||$ $\le C_{2}  \, |||(\bm{u},p_T,{p}_F,p_1,p_2) |||$, satisfying
\begin{align}\label{B_LBound}
  C_{1} \, ||| (\bm{u},p_T,{p}_F,p_1,p_2) |||^{2} \le B(\bm{u},p_T,{p}_F,p_1,p_2; \bm{v},q_T,{q}_F, {q}_1,{q}_2)  \le C_{3}  ||| (\bm{u},p_T,{p}_F,p_1,p_2) |||^{2}
\end{align}
where $C_{1}$, $C_{2}$ and $C_{3}$ depend on the inf-sup constant $C_{D}$ in \eqref{first_inf_sup}, the Korn constant $C_{K}$ and the bounds for $E, \kappa, e_{0}$ and $\kappa_{0}$  stated in Assumption \ref{bounde11}, but not on the physical parameters $\alpha$, $\nu$ or $\tilde{s}_{0}$.
\end{lemma}

\smallskip 

Note that Lemma \ref{main_lemma} implies that the solution $(\bm{u},p_T,{p}_F,p_1,p_2)$ to  \eqref{scm12} satisfies
 $$ C_{1} \,  ||| (\bm{u},p_T,{p}_F,p_1,p_2) |||^{2} \le B(\bm{u},p_T,{p}_F,p_1,p_2; \bm{v},q_T,{q}_F, q_1, q_2) = \bm{f}(\mathbf{v}) + g(q_{F}).$$
Using the definitions of $\bm{f}(\mathbf{v})$ and $g(q_{F})$, one then easily finds that there exists a constant $C_{4}$ (depending only on the Poincar\'e--Freidrichs constant and $e_{0}^{\min}$ and $e_{0}^{\max}$) such that
\begin{align*}
C_{1} \, |||(\bm{u},p_T,{p}_F,p_1,p_2)||| ^{2} \le  & {C}_{4} \left(\tilde{\mu}^{-1/2} ||e_0^{-1/2}\bm{f}||_{\bm{L}^{2}(D)}+||\kappa_0^{-1/2}g||_{L^{2}(D)}\right) \, ||| (\bm{v},q_T,{q}_F, q_1, q_2) |||.
 \end{align*}
We then obtain
 \begin{align}
|||(\bm{u},p_T,{p}_F,p_1,p_2)|||\le {C} \left(\tilde{\mu}^{-1/2} ||e_0^{-1/2}\bm{f}||_{\bm{L}^{2}(D)}+||\kappa_0^{-1/2}g||_{L^{2}(D)} \right),
\label{hadamard}
 \end{align}
 where ${C}:=C_{4}C_{2}/C_{1} >0$ does not depend on $\nu$, $\alpha$ or $\tilde{s}_{0}$.  The proof of Lemma \ref{main_lemma}, which is given below, follows that of \cite[Lemma 2.3]{KPSpre}. However, modifications are needed to deal with the fact that we have five solution fields and we need to work with the weighted norm $||| \cdot |||$  defined in \eqref{normden}. 
 

\smallskip

\begin{proof} \emph{[Lemma \ref{main_lemma}]}. First, we prove the lower bound. From \eqref{big_B}, we have
 \begin{align*}
 B(\bm{u}, \, &p_T,{p}_F,p_1,p_2; \bm{u},-p_T,-{p}_F,p_1,p_2)\\
 &=a(\bm{u},\bm{u})+b(\bm{u},p_T)+b(\bm{u},-p_T)
 -c_1({p}_1,-p_T)+ c_2(p_1,-{p}_F)-c_3(p_2,-p_F)\\
 &\quad+d({p}_F,{p}_F)-c_1(p_T,p_1)+c_2(p_F,p_1)+\tilde{c}_1(p_1,p_1)-c_3(p_F,p_2)+\tilde{c}_2(p_2,p_2),\\
&=a(\bm{u},\bm{u})+ d({p}_F,{p}_F)+\tilde{c}_1(p_1,p_1)+\tilde{c}_2(p_2,p_2)  =:  |\bm{u}|_a^2+ |{p}_F|^2_{d}+|p_1|_{\tilde{c}_1}^2+|p_2|_{\tilde{c}_2}^2.
 \end{align*}
Now, as a consequence of  \eqref{first_inf_sup}, since $p_{T}\in \mathcal{W}$, there exists a $\bm{w}\in \bm{\mathcal{V}}$ such that
\begin{align}\label{infsupcons11}
- b(\bm{w},p_T) \ge C_D  \frac{1}{\tilde{\mu} E_{\max}} ||p_T||^2_{\mathcal{W}},\quad E_{\max}^{1/2}\rbl{\tilde{\mu}}^{1/2}||\nabla \bm{w}||_{\bm{\mathcal{W}}}\le \rbl{\tilde{\mu}}^{-1/2}E_{\max}^{-1/2}||p_T||_{\mathcal{W}}.
\end{align}
Using the following inequalities
\begin{align}
||p_T||^2_{\mathcal{W}}&\ge e_{0}^{\min}||e_0^{-1/2}p_T||^2_{\mathcal{W}},\\
\rbl{\tilde{\mu}}^{1/2}||e_{0}^{1/2}\nabla \bm{w}||_{\bm{\mathcal{W}}}&\le (e_0^{\max})^{1/2}\rbl{\tilde{\mu}}^{1/2}||\nabla \bm{w}||_{\bm{\mathcal{W}}}\le E_{\max}^{1/2}\rbl{\tilde{\mu}}^{1/2}||\nabla \bm{w}||_{\bm{\mathcal{W}}}, \label{A3} \\
\rbl{\tilde{\mu}}^{-1/2}E_{\max}^{-1/2}||p_T||_{\mathcal{W}}&\le \rbl{\tilde{\mu}}^{-1/2}(e_0^{\max})^{-1/2}||p_T||_{\mathcal{W}}\le \rbl{\tilde{\mu}}^{-1/2}||e_0^{-1/2}p_T||_{\mathcal{W}},
\end{align}
in (\ref{infsupcons11}) implies that
\begin{align}\label{b_cond}
&- b(\bm{w},p_T) \ge C_D  \frac{e_0^{\min}}{E_{\max}}\rbl{\tilde{\mu}}^{-1} ||e_0^{-1/2}p_T||^2_{\mathcal{W}},\quad\rbl{\tilde{\mu}}^{1/2}||e_{0}^{1/2}\nabla \bm{w}||_{\bm{\mathcal{W}}}\le  \rbl{\tilde{\mu}}^{-1/2}||e_0^{-1/2}p_T||_{\mathcal{W}}.
\end{align}
Using the chosen $\bm{w}$ in \eqref{big_B} and using \eqref{b_cond} and \eqref{A3}, it follows that\footnote{Here, we used the standard inequality $-ab \ge -\frac{\epsilon a^{2}}{2} - \frac{b^{2}}{2\epsilon}$.}, for any $\epsilon >0$, 
\begin{align*}
 \mathcal{B}(\bm{u},p_T,{p}_F,p_1,p_2; -\bm{w},0,0,0,0)&=\rbl{-b(\bm{w},p_T) -a(\bm{u},\bm{w}) } \\
 &\ge C_D \, \rbl{\tilde{\mu}}^{-1} \frac{e_0^{\min}}{E_{\max}}||e_{0}^{-\frac{1}{2}}p_T||^2_{\mathcal{W}}-|\bm{u}|_{a}|\bm{w}|_{a} \\
 &\ge C_D\,  \frac{e_0^{\min}}{\tilde{\mu}E_{\max}} ||e_0^{-\frac{1}{2}}p_T||^2_{\mathcal{W}} 
             - |\bm{u}|_{a} \, E_{\max}^{1/2} \,  {\tilde{\mu}}^{1/2} ||\nabla \bm{w}||_{\bm{\mathcal{W}}} \\
&\ge C_D \, \frac{e_0^{\min}}{\tilde{\mu}E_{\max}} ||e_0^{-\frac{1}{2}}p_T||^2_{\mathcal{W}}
             -|\bm{u}|_{a} \, {\tilde{\mu}}^{-1/2} ||e_0^{-1/2}p_T||_{\mathcal{W}} \\
&\ge C_D\,  \frac{e_0^{\min}}{\tilde{\mu}E_{\max}}||e_0^{-\frac{1}{2}}p_T||^2_{\mathcal{W}}- \frac{\epsilon} {2}|\bm{u}|_{a}^2
-\frac{1}{2\epsilon\tilde{\mu}} ||e_0^{-\frac{1}{2}}p_T||_{\mathcal{W}}^2,
 \end{align*}

In addition, using \eqref{big_B} again, for any $\epsilon_{1}, \epsilon_{2}, \epsilon_{3} >0$ we have,
 \begin{align*}
  &B(\bm{u},p_T,{p}_F,p_1,p_2; 0,0,0,-p_T,-p_F) \\
   &=c_1(p_T,p_T)-c_2(p_F,p_T)-\tilde{c}_1({p}_1,p_T)+c_3(p_F,p_F)-\tilde{c}_2(p_2,p_F)\\
   &\ge  \tilde{\lambda}^{-1}||p_T||^2_{\mathcal{W}}-\alpha \tilde{\lambda}^{-1}||{p}_F||_{\mathcal{W}}||{p}_T||_{\mathcal{W}}-E_{\max}^{1/2}\tilde{\lambda}^{-1/2}||p_T||_{\mathcal{W}}  |p_1 |_{\tilde{c}_1}  + \tilde{s}_0||p_F||^2_{\mathcal{W}} -E_{\max}^{1/2}\tilde{s}^{1/2}_0||p_F||_{\mathcal{W}}|p_2|_{\tilde{c}_2}, \\
&\ge \frac{1}{\tilde{\lambda}}||p_T||^2_{\mathcal{W}}-\frac{\alpha^2\epsilon_1}{2\tilde{\lambda}}||{p}_F||_{\mathcal{W}}^2- \frac{1}{2\epsilon_1\tilde{\lambda}}||{p}_T||_{\mathcal{W}}^{2} -\frac{E_{\max}}{2\epsilon_2\tilde{\lambda}}||p_T||_{\mathcal{W}}^2-\frac{\epsilon_2}{2}  |p_1|_{\tilde{c}_1}^2 + \tilde{s}_0||p_F||^2_{\mathcal{W}} \\
 &\quad 
 -\frac{E_{\max}\tilde{s}_0}{2\epsilon_3}||p_F||_{\mathcal{W}}^2-\frac{\epsilon_3}{2}  |p_2|_{\tilde{c}_{2}}^{2},\\
&\ge \Big(1-\frac{1}{2\epsilon_1}-\frac{E_{\max}}{2\epsilon_2}\Big)\frac{1}{\tilde{\lambda}}||p_T||^2_{\mathcal{W}}+ \Big(\Big(1-\frac{E_{\max}}{2\epsilon_3}\Big)\tilde{s}_0-\frac{\epsilon_1}{2}\frac{\alpha^2}{\tilde{\lambda}}\Big) ||p_F||^2_{\mathcal{W}}-\frac{\epsilon_2} {2}|{p}_1|_{\tilde{c}_1}^2-\frac{\epsilon_3}{2}|p_2|_{\tilde{c}_2}^2, \\
&\ge \Big(1-\frac{1}{2\epsilon_1}-\frac{E_{\max}}{2\epsilon_2}\Big)\frac{e_0^{\min}}{\tilde{\lambda}}||e_0^{-\frac{1}{2}}p_T||^2_{\mathcal{W}}+ \Big(\Big(1-\frac{E_{\max}}{2\epsilon_3}\Big)\tilde{s}_0-\frac{\epsilon_1}{2}\frac{\alpha^2}{\tilde{\lambda}}\Big) e^{\min}_0||e_0^{-\frac{1}{2}}p_F||^2_{\mathcal{W}} -\frac{\epsilon_2} {2}|{p}_1|_{\tilde{c}_1}^2  -\frac{\epsilon_3}{2}|p_2|_{\tilde{c}_2}^2.
 \end{align*}

Now choose any $\delta >0$ and $\delta' >0$ and consider 
 \begin{align*}
 & B(\bm{u},p_T,{p}_F,p_1,p_2; \bm{u}-\delta \bm{w},-p_T,-p_F,{p}_1-\delta'  p_T, p_2-\delta' p_F)\nonumber\\
 &=B(\bm{u},p_T,{p}_F,p_1,p_2; \bm{u},-p_T,-{p}_F,p_1,p_2)
 +\delta\, B(\bm{u},p_T,{p}_F,p_1,p_2; -\bm{w},0,0,0,0)\nonumber\\
 &\quad+\delta' \,  B(\bm{u},p_T,{p}_F, p_{1}, p_{2} ; 0,0,0, -p_{T},-p_{F}). \nonumber
 \end{align*}
Combining the above lower bounds for the three terms on the right and rearranging gives  
 \begin{align}
 & B(\bm{u},p_T,{p}_F,p_1,p_2; \bm{u}-\delta \bm{w},-p_T,-p_F,{p}_1-\delta'  p_T, p_2-\delta' p_F)\nonumber\\
 & \ge \left(1-\frac{\delta\epsilon}{2}\right)|\bm{u}|_{{a}}^2
 +\Bigg(\frac{\delta}{\tilde{\mu}}\left(\frac{C_De_0^{\min}}{E_{\max}}-\frac{1}{2\epsilon}\right)
 +\frac{\delta'e_0^{\min}}{\tilde{\lambda}}\Big(1-\frac{1}{2\epsilon_1}-\frac{E_{\max}}{2\epsilon_2}\Big)\Bigg)||e_0^{-\frac{1}{2}}p_T||^2_{\mathcal{W}}\nonumber  
 \end{align}
 \begin{align}
 &\quad+|{p}_F|^2_{d}+ {\delta'  \Big(\Big(1-\frac{E_{\max}}{2\epsilon_3}\Big)\tilde{s}_0-\frac{\epsilon_1}{2}\frac{\alpha^2}{\tilde{\lambda}}\Big) e_0^{\min}||e_0^{-\frac{1}{2}}p_F||^2_{\mathcal{W}}}
 +\left(1-\frac{\delta' \epsilon_2}{2}\right)  |{p}_1|_{\tilde{c}_1}^{2}+\left(1-\frac{\delta' \epsilon_3}{2}\right)  |{p}_2|_{\tilde{c}_2}^{2}. \nonumber
 \end{align}
 Next, making the specific choices $\epsilon_1=1, \epsilon_2=\epsilon_3= 2{E_{\max}}$ and
\begin{align}\label{delta_def}
\epsilon= \frac{E_{\max}}{C_De_0^{\min}}, \quad \delta=\frac{1}{\epsilon}=\frac{C_De_0^{\min}}{E_{\max}}, \quad \delta'=\frac{1}{\epsilon_2}=\frac{1}{2E_{\max}},
\end{align}
we have
\begin{align}
 & B(\bm{u},p_T,{p}_F,p_1,p_2; \bm{u}-\delta \bm{w},-p_T,-p_F,{p}_1-\delta'  p_T, p_2-\delta' p_F)\nonumber\\
  &\ge \frac{1}{2}|\bm{u}|_{{a}}^2+ \frac{1}{2} \left(\frac{\delta^2}{\tilde{\mu} }
  +\frac{e_0^{\min}}{4\tilde{\lambda} E_{\max}}\right)||e_0^{-\frac{1}{2}}p_T||^2_{\mathcal{W}}+|{p}_F|^2_{d}+ \frac{e_0^{\min}}{2E_{\max}}\left(\frac{3}{4}\tilde{s}_0-\frac{\alpha^2}{2\tilde{\lambda}}\right) ||e_0^{-\frac{1}{2}}p_F||^2_{\mathcal{W}} + \frac{1}{2}|{p}_1|_{\tilde{c}_1}^2  +\frac{1}{2}|{p}_2|_{\tilde{c}_2}^2. \nonumber
   \end{align} 
 
 To ensure a positive lower bound, we will need to assume that $\frac{3}{4}\tilde{s}_0-\frac{1}{2} \alpha^2 \tilde{\lambda}^{-1} >0$. For this, we need
\begin{align}\label{cond_param}
 0 < \alpha^{2} \tilde{\lambda}^{-1} < \frac{3\tilde{s}_{0}}{2}.
 \end{align}
Let us assume then that $\alpha^2\tilde{\lambda}^{-1}\in(0,\tilde{C}_1\tilde{s}_0]$ with $\tilde{C}_1< 3/2$. This gives
$$\frac{3}{4}\tilde{s}_0-\frac{\alpha^{2}}{2 \tilde{\lambda}} \ge \left(\frac{3}{4} - \frac{\tilde{C}_{1}}{2}\right)\tilde{s}_{0}= : A \tilde{s}_{0} = \frac{2}{5} A \tilde{s}_{0} +\frac{3}{5} A \tilde{s}_{0} \ge  \frac{2}{5} A \tilde{s}_{0}  + \frac{2}{5} A \frac{\alpha^{2}}{\tilde{\lambda}}= \frac{2A}{5}\left(\tilde{s}_{0} +\frac{\alpha^{2}}{\tilde{\lambda}} \right).$$
 If we now define $\tilde{C}:=2A/5$ then we have
 \begin{align}
  B&(\bm{u},p_T,{p}_F,p_1,p_2; \bm{u}-\delta \bm{w},-p_T,-p_F,{p}_1-\delta'  p_T, p_2-\delta' p_F)\nonumber\\
  &\ge \frac{1}{2}|\bm{u}|_{{a}}^2+ \frac{1}{2} \left(\frac{\delta^2}{\tilde{\mu} }
  +\frac{e_0^{\min}}{4\tilde{\lambda} E_{\max}}\right)||e_0^{-\frac{1}{2}}p_T||^2_{\mathcal{W}}+|{p}_F|^2_{d}+ \frac{e_0^{\min}\tilde{C}}{2E_{\max}}\left(\tilde{s}_0+\frac{\alpha^2} {\tilde{\lambda}}\right) ||e_0^{-\frac{1}{2}}p_F||^2_{\mathcal{W}} + \frac{1}{2}|{p}_1|_{\tilde{c}_1}^2  +\frac{1}{2}|{p}_2|_{\tilde{c}_2}^2\nonumber
   \end{align}
  \noindent  and using the bounds \eqref{abound}--\eqref{dubound} gives
  \begin{align*}
   B&(\bm{u},p_T,{p}_F,p_1,p_2; \bm{u}-\delta \bm{w},-p_T,-p_F,{p}_1-\delta'  p_T, p_2-\delta' p_F)\nonumber\\
  &\ge \frac{1}{2} \left( \frac{C_{K} E_{\min}}{e_0^{\max}}\right) \tilde{\mu}\, ||e_0^{\frac{1}{2}}\nabla\bm{u}||_{\bm{\mathcal{W}}}^2 
  + \frac{1}{2 } \left(\frac{\delta^2}{\tilde{\mu} }  +\frac{e_0^{\min}}{4\tilde{\lambda}E_{\max}}\right) ||e_0^{-\frac{1}{2}}p_T||^2_{\mathcal{W}}  +\frac{\kappa_{\min}}{\kappa_0^{\max}}||\kappa_0^{\frac{1}{2}}\nabla p_F||_{\mathcal{W}}^2\nonumber\\
&\quad+  \left(\frac{e_0^{\min}\tilde{C}}{2E_{\max}} \right)(\tilde{s}_0+\alpha^{2} \tilde{\lambda}^{-1}) ||e_0^{-\frac{1}{2}}p_F||^2_{\mathcal{W}}
  + \frac{1}{2}\left(\frac{E_{\min}}{e_{0}^{\max}}\right) \tilde{\lambda}^{-1} ||e_0^{\frac{1}{2}}{p}_1||_{\mathcal{W}}^2 + \frac{1}{2}\left(\frac{E_{\min}}{e_0^{\max}} \right) \tilde{s}_{0} \, ||{e_0^{\frac{1}{2}}p}_2||_{\mathcal{W}}^2,\nonumber\\
  &\ge C_{1}  \left( \tilde{\mu} 
 ||e_0^{\frac{1}{2}}\nabla\bm{u}||_{\bm{\mathcal{W}}}^2 
 +\left(\frac{1}{\tilde{\mu}}+\frac{1}{\tilde{\lambda}} \right)||e_0^{-\frac{1}{2}}p_T||^2_{\mathcal{W}}
 + (\tilde{s}_0+\alpha^2 \tilde{\lambda}^{-1} ) ||e_0^{-\frac{1}{2}}{p}_F||^2_{\mathcal{W}}+||\kappa_0^{\frac{1}{2}}\nabla p_F||_{\mathcal{W}} \right.\nonumber\\
 &\quad\left.+  \tilde{\lambda}^{-1} ||e_0^{\frac{1}{2}}{p}_1||^2_{\mathcal{W}}+\tilde{s}_{0} ||e_0^{\frac{1}{2}}{p}_2||^2_{\mathcal{W}}\right) 
  =: C_{1} \,   ||| (\bm{u}, p_T, {p}_F,p_1,p_2) |||^{2},  \nonumber
 \end{align*}
 where $C_{1} = \frac{1}{2}\min\{ \frac{E_{\min}C_K}{e_0^{\max}}, \frac{C_D^2(e_{0}^{\min})^2}{E_{\max}^2},\frac{e_0^{\min}}{4E_{\max}}, \frac{e_0^{\min}\tilde{C}}{2E_{\max}},\frac{\kappa_{\min}}{\kappa_0^{\max}}\}$ (since $C_{K} \le 1)$. We have shown then that \eqref{B_LBound} holds with $\bm{v}:=\bm{u} - \delta \bm{w},$ $q_T:=-p_T$, $q_F:=-p_F$, ${q}_1 := {p}_1 - \delta' p_T$, and ${q}_2 := {p}_2 - \delta' p_F$.
 
 
 To complete the proof of the lower bound, we note that due to \eqref{b_cond} we have
  \begin{align}
\tilde{\mu} ||e_0^{\frac{1}{2}}\nabla(\bm{u}-\delta\bm{w})||_{\bm{\mathcal{W}}}^2
&\le 2 \tilde{\mu} ||e_0^{\frac{1}{2}}\nabla\bm{u}||_{\bm{\mathcal{W}}}^2
+2 \delta^2 \tilde{\mu}||e_0^{\frac{1}{2}}\nabla\bm{w}||_{\bm{\mathcal{W}}}^2 \le 2 \tilde{\mu} ||e_0^{\frac{1}{2}}\nabla\bm{u}||_{\bm{\mathcal{W}}}^2
+ 2 \delta^2 \tilde{\mu}^{-1} ||e_0^{-\frac{1}{2}}p_T||_{{\mathcal{W}}}^2. \nonumber
 \end{align}
Similarly, using the definition of $\delta'$ gives
 \begin{align*}
\tilde{\lambda}^{-1}  ||e_0^{\frac{1}{2}}({p}_1-\delta' p_T)||_{\mathcal{W}}^2& \le  
2 \tilde{\lambda}^{-1}  \| e_0^{\frac{1}{2}}{p}_1 \|_{\mathcal{W}}^{2}  
+  2  \delta'^{2} \tilde{\lambda}^{-1}  \| e_{0}^{1/2}p_T \|_{\mathcal{W}}^{2},\\
&\le 2 \tilde{\lambda}^{-1}  \| e_0^{\frac{1}{2}}{p}_1 \|_{\mathcal{W}}^{2}  
+   \frac{(e_0^{\max})^2}{2E_{\max}^2}  \, \tilde{\lambda}^{-1} \| e_{0}^{-1/2}p_T \|_{\mathcal{W}}^{2},
 \end{align*}
 and 
\begin{align*}
 \tilde{s}_{0}  ||e_0^{\frac{1}{2}}({p}_2-\delta' p_F)||_{\mathcal{W}}^2&\le 2 \tilde{s}_{0} \| e_0^{\frac{1}{2}}{p}_{2} \|_{\mathcal{W}}^{2}
+   \frac{(e_0^{\max})^2}{2E_{\max}^2} \,  \tilde{s}_{0} \| e_{0}^{-1/2}p_F \|_{\mathcal{W}}^{2}.
 \end{align*}
 Using the definition of $ ||| \cdot |||$ in \eqref{normden} and the fact that $(e_0^{\max})^2/E_{\max}^2 \le 1$ then leads to the upper bound
 \begin{align*} 
 &|||(\bm{u}-\delta\bm{w},-p_T,-p_F,{p}_1-\delta' p_T,p_2-\delta' p_F)|||^2\nonumber\\
 &= \tilde{\mu} ||e_0^{\frac{1}{2}}\nabla(\bm{u}-\delta\bm{w})||^2_{\bm{\mathcal{W}}}
 +\left( \tilde{\mu}^{-1}+ \tilde{\lambda}^{-1}\right)
 ||e_0^{-\frac{1}{2}}p_T||_{\mathcal{W}}^2+ (\tilde{s}_0+\alpha^2 \tilde{\lambda}^{-1})
 ||e_0^{-\frac{1}{2}}p_F||_{\mathcal{W}}^2+
 ||\kappa_0^{\frac{1}{2}}\nabla p_F||_{\mathcal{W}}^2\\
 &\quad+ \tilde{\lambda}^{-1} ||e_0^{\frac{1}{2}}({p}_1-\delta' p_T)||_{\mathcal{W}}^2 +  \tilde{s}_{0} ||e_0^{\frac{1}{2}}({p}_2-\delta' p_F)||_{\mathcal{W}}^2, \nonumber\\
&\le  (2+2 \delta^2) \left( \tilde{\mu} 
 ||e_0^{\frac{1}{2}}\nabla\bm{u}||_{\bm{\mathcal{W}}}^2 
 +\left(\frac{1}{\tilde{\mu}}+\frac{1}{\tilde{\lambda}} \right)||e_0^{-\frac{1}{2}}p_T||^2_{\mathcal{W}}
 + ( \tilde{s}_0+\alpha^2 \tilde{\lambda}^{-1}) ||e_0^{-\frac{1}{2}}{p}_F||^2_{\mathcal{W}}+||\kappa_0^{\frac{1}{2}}\nabla p_F||_{\mathcal{W}} \right. \nonumber\\
 &\quad \left.+   \tilde{\lambda}^{-1} ||e_0^{\frac{1}{2}}{p}_1||^2_{\mathcal{W}} + \tilde{s}_{0} ||e_0^{-\frac{1}{2}}{p}_2||^2_{\mathcal{W}}\right)
  =: C' \, ||| (\bm{u}, p_T, {p}_F,p_1,p_2) |||^{2},
 \end{align*}
 as required. The result now holds with $C_{2}=\sqrt{C'} = \sqrt{2+2\delta^{2}}$ where $\delta$ is defined as in \eqref{delta_def}.
 
To establish the upper bound, we omit the full details for brevity. Clearly,
 \begin{align*}
  B(\bm{u}, \, &p_T,{p}_F,p_1,p_2; \bm{v},q_T, {q}_F, q_1,q_2)\\
 & \le \mid a(\bm{u},\bm{v}) \mid + \mid b(\bm{v},p_T ) \mid +\mid b(\bm{u}, q_T) \mid  + \mid c_1({p}_1, q_T) \mid + \mid c_2(p_1, {q}_F)\mid + \mid c_3(p_2,q_F) \mid + \mid d(p_F, q_F) \mid \\ 
 & +  \mid c_1(p_T,q_1) \mid + \mid c_2(p_F,q_1) \mid + \mid \tilde{c}_1(p_1,q_1) \mid + \mid c_3(p_F,q_2) \mid +\mid \tilde{c}_2(p_2,q_2) \mid.
 \end{align*} 
All of the bilinear forms are bounded. Applying their upper bounds, grouping terms and then applying the Cauchy-Schwarz inequality for sums gives
  \begin{align*}
 B(\bm{u}, \, &p_T,{p}_F,p_1,p_2; \bm{v},q_T, {q}_F, q_1,q_2) \\
 & \le  3 \sqrt{6} \max\left\{\frac{E_{\max}}{e_{0}^{\min}}, \sqrt{d}, \frac{\kappa_{\max}}{\kappa_{0}^{\min}}\right\}  ||| (\bm{u}, p_T, {p}_F,p_1,p_2) |||\,  ||| (\bm{v}, q_T, {q}_F,q_1,q_2) |||.
 \end{align*} 
Using the upper bound for  the second norm, the result holds with $C_{3}:=   3 \sqrt{6}  C_{2} \max\left\{\frac{E_{\max}}{e_{0}^{\min}}, \sqrt{d}, \frac{\kappa_{\max}}{\kappa_{0}^{\min}}\right\}$.
\end{proof}

 
 \smallskip 
 
\begin{remark} Note that \eqref{cond_param} simply says that the physical parameters $\alpha$, $\nu$ and $s_{0}$ need to be chosen in  a compatible way. If $s_{0}=\alpha^{2}/\lambda$ (as assumed in \cite{lee2017parameter}), then $\tilde{s}_{0} = \alpha^{2}\tilde{\lambda}^{-1}$ and  \eqref{cond_param}  is satisfied. If $s_{0}$ is chosen as in \eqref{s0_def}, then $\tilde{s}_{0}$ depends on $\nu, \alpha$ and $\phi$  (but not $E$) and  \eqref{cond_param}  yields a compatibility condition. For example, if $d=2$, we require $\alpha^{2} < 3(\alpha - \phi) \nu$. If  \eqref{cond_param} is not satisfied, then it may still be possible to establish a positive lower bound but with a constant $C_{1}$ that depends on the physical parameters.

\end{remark} 

\section{\rbl{Stochastic Galerkin Mixed Finite Element Method (SG-MFEM)}}\label{disformul}

We now discuss how to construct an SG-MFEM approximation of the solution to (\ref{scm11a}). For this, we first need to choose appropriate finite-dimensional subspaces of $\bm{\mathcal{V}}_{0,\bm{u}}$  and $\mathcal{V}_{0,p}$ to approximate the displacement $\bm{u}$ and fluid pressure $p_{F}$, respectively, and finite-dimensional subspaces of $\mathcal{W}$ to approximate the total pressure $p_{T}$, and auxiliary variables $p_{1}$ and $p_{2}$. As usual, we exploit the fact that $\bm{\mathcal{V}}_{0,\bm{u}} \cong \bm{H}_{0,\bm{u}}^{1}(D) \otimes L_{\pi}^{2}(\Gamma)$ (and similarly for the other spaces) and employ a tensor product construction that combines a mixed finite element method (MFEM) on the spatial domain $D$, and global polynomial approximation on the parameter domain $\Gamma$.


We start by selecting two compatible pairs of finite element spaces associated with a mesh ${\cal T}_{h}$ (with mesh parameter $h$) on the spatial domain $D$. First, we choose a pair $(\bm{V}_{0,\bm{u}}^{h}, W^{h})$ with $\bm{V}_{0,\bm{u}}^{h} \subset \bm{H}_{0,\bm{u}}^{1}(D)$ and $W^{h} \subset L^{2}(D)$ such that the discrete inf--sup condition 
\begin{align}\label{discrete_inf_sup}
\sup_{0\neq\bm{v}\in {\bm{V}_{0,\bm{u}}^{h}}}
\frac{\int_D q \,\nabla\cdot\bm{v} }{ ||\nabla\bm{v}||_{\bm{L}^2(D)} }&\ge \gamma \,  ||q||_{L^2(D)} \quad\forall  q\in W^h
 \end{align}
is satisfied with $\gamma$ uniformly bounded away from zero (i.e., independent of $h$). Next, we choose  a pair $(\widetilde{V}_{0,p}^{h},  \widetilde{W}^h)$ with $\widetilde{V}_{0,p}^{h}  \subset H_{0,p}^{1}(D)$ and $\widetilde{W}^h \subset L^{2}(D)$ such that the discrete inf--sup condition 
\newpage
\begin{align}\label{discrete_inf_sup12}
\sup_{0\neq q\in {\widetilde{W}^{h}}}
\frac{\int_D q \,{p} }{ ||q||_{{L}^2(D)} }&\ge \gamma_L
 ||p||_{L^2(D)} \quad\forall  p\in \widetilde{V}_{0,p}^{h}
 \end{align}
is also satisfied with $\gamma_L$ uniformly bounded away from zero. Clearly (\ref{discrete_inf_sup12}) always holds with  $\gamma_L=1$ if $\widetilde{V}_{0,p}^{h} \subset \widetilde{W}^{h}$. We can ensure this is true by choosing the same $H^{1}(D)$-conforming finite element space for both $\widetilde{V}_{0,p}^{h}$ and $\widetilde{W}^{h}$  (but removing basis functions associated with nodes on $\partial D_{p}$ to define $\widetilde{V}_{0,p}^{h}$).
%
%

Now, let $Q_{1}$ (resp., $Q_{2}$) denote the usual finite element space associated with continuous piecewise bilinear (resp., biquadratic) approximation and let $\bm{Q}_{1}$ (resp., $\bm{Q}_{2}$) denote the analogous vector-valued space, with components in $Q_{1}$ (resp.,  $Q_{2}$). Some possible combinations $\bm{V}_{0,\bm{u}}^{h}-W^{h}-\widetilde{V}_{0,p}^{h}-W^{h}-\widetilde{W}^{h}$  of finite element spaces  for the spatial approximation of $\bm{u}$, $p_{T}$, $p_{F}$, $p_{1}$ and $p_{2}$ (in that order), are as follows. \medskip 

\begin{itemize}
\item {$\bm{Q}_2-Q_1-Q_1-Q_1-Q_1$:} This method uses the standard Taylor--Hood $\bm{Q}_2-Q_{1}$ element  for the spatial approximation of $\bm{u}$  and $p_T$.  It is well-known (e.g. see \cite{HDA})  that (\ref{discrete_inf_sup}) is satisfied for this pair. Since $p_F$ and $p_2$ are both approximated in $Q_1$ (up to boundary conditions), \eqref{discrete_inf_sup12} is also trivially satisfied. However, a priori error analysis for $\bm{Q}_{2}-Q_{1}-Q_{1}$ approximation for the deterministic three-field model in \cite{oyarzua2016locking} suggests that if $\bm{u}, p_{F}$ and $p_{T}$ have sufficient spatial regularity, one can only expect $O(h)$ convergence in an energy-type norm. 

\smallskip 

\item {$\bm{Q}_2-Q_1-Q_2-Q_1-Q_2$:} This method uses $Q_{2}$ approximation for $p_F$ and $p_2$ as well as for the displacement.  Both inf-sup conditions are satisfied and the a priori results in \cite{oyarzua2016locking} suggest that one can now expect $O(h^{2})$ convergence in an energy norm if the solution has sufficient spatial regularity.
  

 \smallskip 
 
 \item $\bm{Q}_{k+1}-Q_k-Q_k-Q_k-Q_k$ or $\bm{Q}_{k+1}-Q_k-Q_{k+1}-Q_k-Q_{k+1}$ for $k \ge 2$:
 For these higher order methods, both inf-sup conditions are satisfied.
 \end{itemize}

\medskip 


Next, we describe the parametric approximation. Recall that $M:=M_{1}+M_{2}$ is the total number of stochastic parameters in \eqref{E_def}--\eqref{k_def}. For each $j=1,\ldots, M,$ we find a set $\{ \psi_{i}(Y_{j}), i=0,1, \ldots \}$ of univariate polynomials on $\Gamma_{j}$ (where $\psi_{i}$ has degree $i$) that are orthonormal in the $L_{\pi_{j}}^{2}(\Gamma_{j})$-sense.  For example, if $Y_{j}$ is the image of $\xi_{j} \sim U(-1,1)$, it is natural to choose $\pi_{j}$ to be the associated probability measure and the required polynomials are Legendre polynomials. Given $M$ appropriate sets of univariate polynomials, and a set of multi-indices $ \Lambda \subset \mathbb{N}_{0}^{M}$, a set of \emph{multivariate} polynomials on $\Gamma$ can then be constructed as follows 
\begin{align}
S_{\Lambda}: =\textrm{span}\left\{ \psi_{\boldsymbol{\alpha}}(\bm{Y})
 = \prod_{i=1}^{M} \psi_{\alpha_{i}}(Y_{i}), \quad \boldsymbol{\alpha} 
 \in \Lambda \right\} \subset L_{\pi}^{2}(\Gamma).
\end{align}
If $\pi$ is a product measure, the basis functions for $S_{\Lambda}$ are by construction \emph{orthonormal} in the $L_{\pi}^{2}(\Gamma)$ sense.  Working with an orthonormal basis is advantageous because the associated matrix components of the Galerkin system are sparse (see Section \ref{matrices}). 



Once we have chosen four suitable finite element spaces $\bm{V}_{0,\bm{u}}^{h}$, $W^{h}$, $\widetilde{V}_{0,p}^{h}$, $\widetilde{W}^{h}$ and an appropriate set of polynomials $S_{\Lambda}$ (equivalently, a set of multi-indices $\Lambda$), we can define the SG-MFEM spaces 
$$\bm{V}_{h,\Lambda}^{0}:=\bm{V}_{0,\bm{u}}^{h} \otimes S_{\Lambda}, \quad W_{h,\Lambda}:=W^{h} \otimes S_{\Lambda}, \quad \widetilde{V}_{h,\Lambda}^{0}:= \widetilde{V}_{0,p}^{h}\otimes S_{\Lambda}, \quad \widetilde{W}^{h,\Lambda}:= \widetilde{W}^{h}\otimes S_{\Lambda},$$
and solve the discrete weak problem: find
 $(\bm{u}^{h,\Lambda}, p^{h,\Lambda}_{T},$ ${p}^{h, \Lambda}_F, p_1^{h,\lambda},p_2^{h,\Lambda})
  \in \bm{V}_{h, \Lambda}^{0} \times W_{h,\Lambda} \times \widetilde{V}_{h,\Lambda}^{0}\times W_{h,\Lambda}\times \widetilde{W}_{h,\Lambda}$ such that
\begin{subequations} \label{dscm11a}
\begin{align}
a(\bm{u}^{h,\Lambda},\bm{v})+b(\bm{v},p_T^{h,\Lambda})&=f(\bm{v}) \quad \quad  \, \, \, \,  \forall \bm{v}\in \bm{{V}}_{h,\Lambda}^{0},\\
b(\bm{u}^{h,\Lambda},q_T)-c_1({p}_1^{h,\Lambda},q_T)&=0 \quad \quad \quad \quad \, \forall q_T\in {W}_{h,\Lambda},\\
c_2(p_1^{h,\Lambda},{q}_F)-c_3(p_2^{h,\Lambda},q_F)-d(p_F^{h,\Lambda}, q_F)&=g(q_F) \, \, \quad \quad \forall  {q}_F\in \widetilde{V}_{h,\Lambda}^{0},\\
-c_1(p_T^{h,\Lambda},q_1)+c_2(p_F^{h,\Lambda},q_1)+\tilde{c}_1(p_1^{h,\Lambda},q_1)&=0\,\,\quad\quad \qquad \forall q_1\in {W}_{h,\Lambda}, \\
-c_3(p_F^{h,\Lambda},q_2)+\tilde{c}_2(p_2^{h,\Lambda},q_2)&=0\,\,\quad\quad \qquad \forall q_2\in \widetilde{W}_{h,\Lambda}.
\end{align}
\end{subequations}
If Assumption \ref{Assump2} is satisfied and if $\bm{V}_{0,\bm{u}}^{h}$, $W^{h}$, $\widetilde{V}_{0,p}^{h}$, $\widetilde{W}^{h}$ are chosen so that the discrete inf-sup conditions (\ref{discrete_inf_sup}) and  (\ref{discrete_inf_sup12}) are satisfied, then well-posedness of \eqref{dscm11a} can be established in a similar way as for \eqref{scm12}.

\subsection{\rbl{SG-MFEM Linear System}}\label{matrices}

We now describe the structure of the discrete linear system associated with \eqref{dscm11a}. We will assume that $\widetilde{V}_{0,p}^{h}$ and $\widetilde{W}_{h}$ are both $Q_{1}$ spaces, or are both $Q_{2}$ spaces so that (\ref{discrete_inf_sup12}) holds. To simplify the description, we will further assume that $\widetilde{W}_{h}=W_{h}$. 

 First, let $n_{y}:=\textrm{dim}(S_{\Lambda})$ and define the symmetric matrices $G_{0}, G_{k} \in \mathbb{R}^{n_{y} \times n_{y}}$ for $k=1,\ldots, M$ by
  \begin{align*}
[G_{0}]_{\boldsymbol{\alpha},\boldsymbol{\beta}} := \int_{\Gamma}\psi_{\boldsymbol{\alpha}}(\bm{Y}) 
\psi_{\boldsymbol{\beta}}(\bm{Y}) \, d\pi(\bm{Y}), \qquad   [G_{k}]_{\boldsymbol{\alpha},\boldsymbol{\beta}} 
:= \int_{\Gamma}  \, Y_{k} \,  \psi_{\boldsymbol{\alpha}}(\bm{Y}) \psi_{\boldsymbol{\beta}}(\bm{Y}) \, d\pi(\bm{Y}), \quad \boldsymbol{\alpha}, \boldsymbol{\beta} \in \Lambda.
  \end{align*}
Since the basis functions for $S_{\Lambda}$ are orthonormal we have  $G_{0}=I$. In addition (due to the three-term recurrence of the underlying orthogonal polynomials  \cite{powell_elman}), $G_{k}$ has at most two \rbl{nonzero} entries per row.  

Now, let $V_{0,\bm{u}}^{h}:= \textrm{span}\{\phi_{1}(\bm{x}), \ldots,  \phi_{n_{u}}(\bm{x})\} \subset H_{0,\bm{u}}^{1}(D)$ denote the usual (scalar-valued) $Q_{1}$ or $Q_{2}$ finite element space associated with mesh nodes that \emph{do not} lie on $\partial D_{\bm{u}}$.  When $D \subset \mathbb{R}^{2}$, a basis for $\bm{V}_{0,\bm{u}}^{h} \subset \bm{H}_{0,\bm{u}}^{1}(D)$, with components in $V_{0,\bm{u}}^{h}$, is then given by
$$\left\{ \left(\begin{array}{c} \phi_{i}(\bm{x})  \\ 0 \end {array} \right), \left(\begin{array}{c} 0 \\ \phi_{i}(\bm{x}) \end {array} \right), \quad  i=1,\ldots, n_{u}\right\}.$$
Using this basis, we define the matrix $A_{11}^{k}  \in \mathbb{R}^{n_{u} \times n_{u}}$ by
 \begin{eqnarray*}
 [A_{11}^{k}]_{i,\ell}  : =   & \int_{D} e_{k}(\bm{x})  \,  \bm{\epsilon}\left(
 \begin{array}{c} \phi_i(\bm{x}) \\ 0 \end{array} \right) : \bm{\epsilon}\left(\begin{array}{c} \phi_\ell(\bm{x}) \\ 0 
 \end{array} \right) \, d\bm{x},   \quad\quad i,\ell=1, \dots, n_{u},
 \end{eqnarray*}
  for $k=0,1, \ldots, M_{1}$ and the matrix $A_{21}^{k} \in \mathbb{R}^{n_{u} \times n_{u}}$ by
 \begin{align*}
[A_{21}^{k}]_{i,\ell} := \int_{D} e_{k}(\bm{x})  \,  \bm{\epsilon}\left(
\begin{array}{c}0 \\ \phi_i(\bm{x})  \end{array} \right) : \bm{\epsilon}\left(\begin{array}{c} \phi_\ell(\bm{x}) \\ 0 
\end{array} \right) \, d\bm{x}, \quad\quad i,\ell=1, \dots, n_{u},
  \end{align*}
and similarly for $A_{12}^{k}, A_{22}^{k} \in \mathbb{R}^{n_{u} \times n_{u}}$, for $k=0,1, \ldots, M_{1}$.

Next, let   $\widetilde{V}_{0,h}^{p}= \textrm{span}\left\{  \varphi_{1}(\mathbf{x}), \ldots,   \varphi_{n_{0}}(\mathbf{x}) \right\} \subset H_{0,p}^{1}(D)$ and let 
$$\widetilde{W}_{h} = W_{h}=\textrm{span}\left\{  \varphi_{1}(\mathbf{x}), \ldots,   \varphi_{n_{0}}(\mathbf{x}), \varphi_{n_{0}+1}(\mathbf{x}), \ldots, \varphi_{n_{p}}(\mathbf{x}) \right\}  \subset H^{1}(D)$$
(so $ \varphi_{n_{0}+1}, \ldots, \varphi_{n_{p}}$ are basis functions associated with nodes on $\partial D_{p}$), and define  $B_{1}, B_{2} \in \mathbb{R}^{n_{p} \times n_{u}}$ by  \begin{align*}
  [B_{1}]_{r,\ell} = - \int_{D}   \varphi_{r}(\mathbf{x}) \,  \frac{\partial \phi_{\ell}(\mathbf{x})}{\partial x_{1}} d \bm{x} , \quad   [B_{2}]_{r,\ell} = - \int_{D}   \varphi_{r}(\mathbf{x}) \,  \frac{\partial \phi_{\ell}(\mathbf{x})}{\partial x_{2}} d \bm{x},
  \end{align*}
  for $r=1,\ldots, n_{p},  \, \ell=1, \ldots, n_{u}$. We also define the mass matrix $C \in \mathbb{R}^{n_{p} \times n_{p}}$ associated with $\widetilde{W}_{h}=W_{h}$ by
    \begin{align*}
  [C]_{r,s} = \int_{D} \varphi_{r}(\mathbf{x}) \, \varphi_{s}(\mathbf{x}) \, d\mathbf{x} , \qquad r,s=1,\ldots, n_{p},
  \end{align*}
  and for $k=0,1,\ldots,M_1$, we define the \emph{weighted} mass matrices $\tilde{C}_{k} \in \mathbb{R}^{n_{p} \times n_{p}}$ by
  \begin{align*}
  [\tilde{C}_{k}]_{r,s} =  \int_{D} e_{k}(\mathbf{x}) \varphi_{r}(\mathbf{x})\, \varphi_{s}(\mathbf{x}) \, d\mathbf{x},
   \qquad r,s=1,\ldots, n_{p}.
  \end{align*}
  The rectangular matrix $C_{b} \in \mathbb{R}^{n_{0} \times n_{p}}$ associated with $\widetilde{V}_{0,p}^{h}$ and $\widetilde{W}_{h}$ is formed by taking the matrix $C$ and simply removing the rows associated with nodes on $\partial D_{p}$.  For each $k=0,1,\ldots,M_2$, we also define the \emph{weighted} stiffness matrices $D_k\in \mathbb{R}^{n_{0} \times n_{0}}$ associated with $\widetilde{V}_{0,p}^{h}$ by 
 \begin{align*}
  [D_{k}]_{r,s} =  \int_{D} \kappa_{k}(\mathbf{x}) \nabla\varphi_{r}(\mathbf{x})\,\cdot\, \nabla\varphi_{s}(\mathbf{x}) \, d\mathbf{x},
   \qquad r,s=1,\ldots, n_{0}.
  \end{align*}
  
The vector $\bm{g}_{0} \in \mathbb{R}^{n_{y}}$ is defined to be the first column of $G_{0}=I$.  Writing $\bm{f}(\bm{x})=(f_{1}(\bm{x}), f_{2}(\bm{x}))^{\top}$, we also define the vectors $\bm{f}_{1}, \bm{f}_{2} \in \mathbb{R}^{n_{u}}$ by
 \begin{align*}
 [\bm{f}_{1}]_{\ell} = \int_{D} f_{1}(\bm{x}) \phi_{\ell}(\bm{x}) \, d \bm{x}, 
 \qquad  [\bm{f}_{2}]_{\ell} = \int_{D} f_{2}(\bm{x}) \phi_{\ell}(\bm{x}) \, d \bm{x}, \qquad \ell=1,\ldots, n_{u},
 \end{align*}
 and finally, we define $\bm{g}\in \mathbb{R}^{n_{0}}$ by
 \begin{align*}
 [\bm{g}]_r=\int_D g(\bm{x})\varphi_{r} (\bm{x}) \, d\bm{x}, \qquad r=1,\ldots, n_{0}.
 \end{align*}
 




Now, permuting the variables in \eqref{dscm11a} so that they appear in the order $\bm{u}^{h,\Lambda}$, ${p}_{1}^{h, \Lambda}$, $p_{2}^{h,\Lambda}$, $p_{F}^{h,\Lambda}$ and $p_{T}^{h, \Lambda}$ leads to a system of $(2n_{u}+3n_{p}+n_{0})n_{y}$ equations with the saddlepoint structure
\begin{align}
\label{saddlepoint}
\left(\begin{array}{cc} 
 \mathcal{A} & \;\;\;\mathcal{B}^{\top} \\ 
 \mathcal{B}  & -\mathcal{C}  \end{array} \right)\left(\begin{array}{c} \mathbf{v} \\ \mathbf{p} \end{array} \right) = \left(\begin{array}{c} \mathbf{b} \\ \mathbf{c} \end{array} \right).
\end{align} 
Here, the solution vector has the block structure
$$ \mathbf{v}= \left(\begin{array}{c} \mathbf{u}_{1} \\ \mathbf{u}_{2} \\ {\mathbf{p}_1}\\ \mathbf{p}_2 \end{array} \right), 
\qquad
\mathbf{p}=\left(\begin{array}{c}\mathbf{p}_F \\ \mathbf{p}_T \end{array}\right),$$
where $\mathbf{u}_{1,2}\in \mathbb{R}^{n_{u}n_{y}}$, ${\mathbf{p}}_1, \mathbf{p}_2, \mathbf{p}_{T} \in \mathbb{R}^{n_{p}n_{y}}$ and $\mathbf{p}_F \in \mathbb{R}^{n_{0}n_{y}}$ contain the degrees of freedom associated with the distinct physical variables,  and on the right-hand side,
$$ \mathbf{b}= \left(\begin{array}{c} \mathbf{g}_{0} \otimes \mathbf{f}_{1} \\ \mathbf{g}_{0} \otimes \mathbf{f}_{2} \\ \mathbf{0} \\ \mathbf{0} \end{array} \right) \in \mathbb{R}^{2(n_{u} + n_{p}) n_{y}}, \qquad 
\mathbf{c}=\left(\begin{array}{c} \mathbf{g}_{0} \otimes \mathbf{g}\\ \mathbf{0}\end{array}\right) \in \mathbb{R}^{(n_{0} + n_{p}) n_{y}}. $$

To obtain \eqref{saddlepoint}, the degrees of freedom must first be grouped by physical variable. For each physical variable, the spatial degrees of freedom are then grouped for the same parametric basis function $\psi_{\boldsymbol{\alpha}(j)}$.  For example, the vector associated with the fluid pressure $p_{F}^{h, \Lambda}$ has the form $\mathbf{p}_{F}^{\top} = (\mathbf{p}_{F,1}^{\top}, \mathbf{p}_{F,2}^{\top}, \ldots, \mathbf{p}_{F,n_{y}}^{\top})$ where $\mathbf{p}_{F,j} \in \mathbb{R}^{n_{0}}$ for $j=1,\ldots, n_{y}$. With the assumed ordering of the degrees of freedom, the blocks of the coefficient matrix in \eqref{saddlepoint} are then given by
\smallskip 
 \begin{align}\label{defexacta}
\mathcal{A} := & {\small \left(\begin{array}{cc|c|c} 
\tilde{\mu}  \sum\limits_{k=0}^{M_1}  G_{k} \otimes A_{11}^{k}   & \tilde{\mu} \sum\limits_{k=0}^{M_1} G_{k} \otimes A_{21}^{k} & \mathbf{0}&\mathbf{0} \\
& &  \\
\tilde{\mu}  \sum\limits_{k=0}^{M_1} G_{k} \otimes A_{12}^{k}  & \tilde{\mu}  \sum\limits_{k=0}^{M_1}  G_{k} \otimes A_{22}^{k}   & \mathbf{0}&\mathbf{0}   \\ \hline 
 \mathbf{0} & \mathbf{0} &   {\tilde{\lambda}}^{-1}\sum\limits_{k=0}^{M_1} G_{k} \otimes \tilde{C}_{k}&\mathbf{0} \\ \hline
  \mathbf{0} & \mathbf{0} & \mathbf{0} & {\tilde{s}}_{0} \sum\limits_{k=0}^{M_1} {G_{k} \otimes \tilde{C}_{k}}
  \end{array} \right)},
\end{align}
and
  \begin{align}\label{defexactb}
\mathcal{B} :=  {\small  \left(\begin{array}{cc|c|c} 
\mathbf{0} & \mathbf{0} & \alpha\tilde{\lambda}^{-1} I \otimes C_{b} & \tilde{s}_{0} I \otimes C_{b}\\
& &\\
 I \otimes B_{1} & I \otimes B_{2} &   -\tilde{\lambda}^{-1} I \otimes C & \mathbf{0}
 \end{array} \right)}, \qquad 
\mathcal{C} :=  {\small \left(\begin{array}{cc}
\sum\limits_{k=0}^{M_2}\tilde{G}_k\otimes D_k  & \mathbf{0}\\ \smallskip
\mathbf{0} & \mathbf{0}
 \end{array} \right)},
 \end{align}
where $\tilde{G}_0=G_0=I$ (the $n_y \times n_y$ identity matrix) and $\tilde{G}_{k}=G_{k+M_1}$ for $k=1,\ldots,M_2$. 


One can write the discrete system in various ways. For instance, by \emph{reordering} the degrees of freedom, so that all of the spatial degrees of freedom for all of the physical variables are grouped for the same parametric basis function $\psi_{\boldsymbol{\alpha}(j)}$, we can also write the linear system in so-called Kronecker form as
\begin{align}\label{alt_system}
\left(G_{0} \otimes {\cal K}_{0} + \sum_{k=1}^{M_{1}} G_{k} \otimes {\cal K}_{k} + \sum_{k=1}^{M_{2}} \tilde{G}_{k} \otimes \tilde{{\cal K}}_{k}\right)\mathbf{x} = \mathbf{z}
\end{align}
where the solution vector has the form $\mathbf{x} = (\mathbf{x}_{1}^{\top}, \ldots, \mathbf{x}_{n_{y}}^{\top})^{\top}$, with
$$ \mathbf{x}_{j}^{\top}: =(\mathbf{u}_{1,j}^{\top}, \mathbf{u}_{2,j}^{\top}, \mathbf{p}_{1,j}^{\top}, \mathbf{p}_{2,j}^{\top}, \mathbf{p}_{F,j}^{\top}, \mathbf{p}_{T,j}^{\top} )^{\top}, \qquad j=1,\ldots, n_{y}.$$
In \eqref{alt_system} we then have 
\begin{align}\label{bigK0}
\mathcal{K}_{0} := & {\small \left(\begin{array}{cccc|cc} 
\tilde{\mu} A_{11}^{0}   & \tilde{\mu} A_{21}^{0} & \mathbf{0}&\mathbf{0}  & \mathbf{0}  & B_{1}^{\top} \\
\tilde{\mu} A_{12}^{0}  & \tilde{\mu}  A_{22}^{0}   & \mathbf{0}&\mathbf{0} &  \mathbf{0}  & B_{2}^{\top}  \\
 \mathbf{0} & \mathbf{0} &   {\tilde{\lambda}}^{-1} \tilde{C}_{0}&\mathbf{0} & \alpha\tilde{\lambda}^{-1} C_{b}^{\top}  &  -\tilde{\lambda}^{-1} C  \\ 
  \mathbf{0} & \mathbf{0} & \mathbf{0} & {\tilde{s}}_{0} {\tilde{C}_{0}} &   \tilde{s}_{0} C_{b}^{\top} &  \mathbf{0}  \\ \hline 
  \mathbf{0} & \mathbf{0} & \alpha\tilde{\lambda}^{-1} C_{b} & \tilde{s}_{0} C_{b} & D_0  & \mathbf{0}\\
 B_{1} &  B_{2} &   -\tilde{\lambda}^{-1} C & \mathbf{0} &\mathbf{0} & \mathbf{0}
  \end{array} \right)}, \qquad \mathbf{z} =  \mathbf{g}_{0} \otimes \left(\begin{array}{c} \mathbf{f}_{1} \\ \mathbf{f}_{2} \\ \mathbf{0} \\ \mathbf{0} \\  \mathbf{g} \\ \mathbf{0} \end{array} \right)
\end{align}
and for $k=1,\ldots, M_{1}$ and $\ell= 1, \ldots, M_{2},$ we have
\begin{align}\label{bigKk}
\mathcal{K}_{k} := & {\small \left(\begin{array}{cccc|cc} 
\tilde{\mu}  A_{11}^{k}   & \tilde{\mu}  A_{21}^{k} & \mathbf{0}&\mathbf{0}  & \mathbf{0}  & \mathbf{0} \\
\tilde{\mu} A_{12}^{k}  & \tilde{\mu} A_{22}^{k}   & \mathbf{0}&\mathbf{0} &  \mathbf{0}  &  \mathbf{0} \\ 
 \mathbf{0} & \mathbf{0} &   {\tilde{\lambda}}^{-1} \tilde{C}_{k}&\mathbf{0} & \mathbf{0}  & \mathbf{0}  \\  
  \mathbf{0} & \mathbf{0} & \mathbf{0} & {\tilde{s}}_{0} {\tilde{C}_{k}} & \mathbf{0} &  \mathbf{0}  \\ \hline
   \mathbf{0} & \mathbf{0} & \mathbf{0} &\mathbf{0} & \mathbf{0}  & \mathbf{0}\\ 
\mathbf{0} & \mathbf{0}&  \mathbf{0} & \mathbf{0} &\mathbf{0} & \mathbf{0}
  \end{array} \right)}, \qquad   \tilde{\mathcal{K}}_{\ell} :=  {\small \left(\begin{array}{cccc|cc} 
\mathbf{0}   & \mathbf{0}   & \mathbf{0}&\mathbf{0}  & \mathbf{0}  & \mathbf{0} \\
\mathbf{0}   & \mathbf{0}     & \mathbf{0}&\mathbf{0} &  \mathbf{0}  &  \mathbf{0} \\ 
 \mathbf{0} & \mathbf{0} &  \mathbf{0}   &\mathbf{0} & \mathbf{0}  & \mathbf{0}  \\ 
  \mathbf{0} & \mathbf{0} & \mathbf{0} &\mathbf{0}  & \mathbf{0} &  \mathbf{0}  \\ \hline
  \mathbf{0} & \mathbf{0} & \mathbf{0} &\mathbf{0} & D_{\ell}  & \mathbf{0}\\
\mathbf{0} & \mathbf{0}&  \mathbf{0} & \mathbf{0} &\mathbf{0} & \mathbf{0}
  \end{array} \right)}.
\end{align}
Using the definition of the Kronecker product, one can also rewrite \eqref{alt_system} as a matrix equation \cite{Valeria},
\begin{align}\label{mat_eq}
K_{0} X + \sum_{k=1}^{M_{1}} {\cal K}_{k}X G_{k} + \sum_{\ell=1}^{M_{2}} \tilde{\cal K}_{\ell} X \tilde{G}_{\ell}= Z
\end{align}
where $X$ is the $(2n_{u} + 3 n_{p} + n_{0}) \times n_{y}$ solution matrix obtained by reshaping the vector $\mathbf{x}$ (so the $j$th column of $X$ is $\mathbf{x}_{j}$) and similarly for the right-hand side matrix $Z$.

Due to their extremely large size, it is usually infeasible to form the Kronecker products that appear in the above equations. However, the coefficient matrices in \eqref{saddlepoint} and \eqref{alt_system} are \emph{symmetric} and \emph{indefinite} so one can in principle solve these systems iteratively using the minimal residual method (MINRES) see~\cite[Chapter 4]{HDA} if (i) enough memory is available to store vectors of length $n_{x}n_{y}$, where $n_{x}:=2n_{u}+3n_{p}+n_{0}$ and (ii) matrix-vector products can be efficiently computed. However, preconditioning is essential due to the inherent ill-conditioning that stems from the discretisation and physical parameters. We discuss this next. When memory is exhausted and  vectors of length $n_{x}n_{y}$ \emph{cannot} be stored, low-rank or reduced-basis solvers (see \cite{Valeria}, \cite{MultiRB}) that operate on \eqref{mat_eq} may need to be explored instead of standard Krylov methods. 

\section{Preconditioning} \label{precond-sec}

We will assume the degrees of freedom are ordered so that the linear system has the form \eqref{saddlepoint}. A natural starting point (see \cite{ernst2009efficient}, \cite{klawonn1998optimal}, \cite{mardal2011preconditioning}, \cite{ss11}) is then to consider block-diagonal preconditioners $P$ with
  \begin{align}\label{Prec_diag}
P = \left(\begin{array}{cc} P_{\mathcal{A}} & 0 \\ 0 & P_{\mathcal{S}}
 \end{array} \right),
 \end{align}
where $P_{\mathcal{A}}$ and $P_{\mathcal{S}}$ approximate $\mathcal{A}$ and $\mathcal{S} = \mathcal{B}\mathcal{A}^{-1} \mathcal{B}^{\top}+\mathcal{C}$ (the Schur complement), respectively. Ideally, we want to choose matrices that are spectrally equivalent to $\mathcal{A}$ and  $\mathcal{S}$,  so that we obtain spectral bounds for the preconditioned system that are independent of the SG-MFEM discretisation parameters \emph{and} the physical parameters $\nu$, $\alpha$ and $\tilde{s}_{0}$.  We can achieve this by following the operator approach to constructing preconditioners described in \cite{mardal2011preconditioning} and \cite{lee2017parameter}. Briefly, the main idea is as follows. 


In Lemma \ref{main_lemma} we established that the solution to \eqref{scm12} is bounded with respect to the norm in \eqref{normden}.  This  result can be used to show (see \cite{mardal2011preconditioning}) that the operator $\mathit{B}$ associated with \eqref{scm12} is a bounded linear map from $\bm{\mathcal X}$ (equipped with $||| \cdot |||$) to its dual space $\bm{\mathcal X}^{*}$ and has a well-defined inverse  $\mathit{B}^{-1}: \bm{\mathcal X}^{*} \to \bm{\mathcal X}$. Moreover,  the operator norms satisfy $\| \mathit{B}\|_{\mathcal{L}(\bm{\mathcal X}, \bm{\mathcal X}^{*})}\le C_{3}$ and $\| \mathit{B} ^{-1}\|_{\mathcal{L}(\bm{\mathcal X}, \bm{\mathcal X}^{*})} \le C_{2}/C_{1}$ and hence the condition number $\kappa(\mathit{B}) \le C_{3}C_{2}/C_{1}$ is bounded independently of  $\nu, \alpha$ and $\tilde{s}_{0}$.  An optimal and  `parameter-robust' preconditioning operator is given by the Riesz map ${\mathcal{R}}: \bm{\mathcal X}^{*} \to \bm{\mathcal X}$ since it can be shown that $\kappa( \mathcal{R}\mathit{B}) \le C_{3}C_{2}/C_{1}$. An analogous result to Lemma \ref{main_lemma} holds for \eqref{dscm11a}, except that the constants $C_{1}, C_{2}, C_{3}$ appearing in the bounds depend on the \emph{discrete} inf-sup constants rather than $C_{D}$. If inf-sup stable FEM spaces are chosen, then the operator associated with \eqref{dscm11a} is a bounded linear map from $\bm{\mathcal X}_{h,\Lambda}$ to  $\bm{\mathcal X}_{h,\Lambda}^{*}$ where  $\bm{\mathcal X}_{h,\Lambda}: = \bm{V}_{h, \Lambda}^{0} \times W_{h,\Lambda} \times \widetilde{V}_{h,\Lambda}^{0}\times W_{h,\Lambda}\times \widetilde{W}_{h,\Lambda}$. The associated operator norm, that of the inverse operator, and the condition number are all then bounded independently of the parameters $\nu, \alpha$ and $\tilde{s}_{0}$ as well as the mesh parameter $h$ and other discretisation parameters  involved in the definition of $S_{\Lambda}$. A `parameter-robust' preconditioner for the finite-dimensional problem is then given by the Riesz map from $\bm{\mathcal X}_{h,\Lambda}^{*}$ to $\bm{\mathcal X}_{h,\Lambda}$. To represent this in matrix form, we choose $P_{\mathcal{A}}$ and $P_{\mathcal{S}}$ in \eqref{Prec_diag} so that, for any $(\bm{v}, q_{1}, q_{2}, q_{F}, q_{T}) \in \bm{\mathcal X}_{h,\Lambda},$
\begin{align} \label{disc_norm_rep}
\bm{x}^{\top} P \bm{x} =  ||| (\bm{v}, q_{1}, q_{2}, q_{F}, q_{T}) |||^{2} 
\end{align} 
where $\bm{x}=(\bm{w}^{\top}, \bm{q}_{1}^{\top}, \bm{q}_{2}^{\top}, \bm{q}_{F}^{\top}, \bm{q}_{T}^{\top}) \in \mathbb{R}^{n_{x}n_{y}}$ is the vector of coefficients associated with  $(\bm{v}, q_{1}, q_{2}, q_{F}, q_{T})$ when the components are expanded in the chosen bases.






 \subsection{Approximation of $\mathcal{A}$}
First, define $\mathbb{A}:=2(A_{11}^{0}+A_{22}^0)/3$ and consider the block-diagonal matrix
 \begin{align}\label{meanpre}
P_{\mathcal{A}} := & {\small  \left(\begin{array}{cc|c|c} 
\tilde{\mu}  \,  I \otimes \mathbb{A}  & \bm{0}  & \mathbf{0} & \mathbf{0}\\
\bm{0}   & \tilde{\mu} \, I \otimes \mathbb{A}   & \mathbf{0}  & \mathbf{0}   \\ \hline  
& & & \\
 \mathbf{0} & \mathbf{0} & \tilde{\lambda}^{-1} I \otimes  \tilde{C}_{0}  & \mathbf{0} \\  \hline 
 & & & \\
  \mathbf{0}& \mathbf{0} & \mathbf{0} & \tilde{s}_{0} I \otimes \tilde{C}_0 
 \end{array} \right)}    =:  {\small  \left(\begin{array}{ccc} 
P_{\mathcal{A},1}& \bm{0} & \bm{0}     \\  
 \mathbf{0} &  P_{\mathcal{A},2} & \mathbf{0}  \\ 
 \mathbf{0}&  \mathbf{0} & P_{\mathcal{A},3} 
 \end{array} \right)}.
 \end{align}
 For any $\bm{v} \in \bm{V}_{h,\Lambda}^{0}$ and $q_{1}, q_{2} \in W_{h,\Lambda},$ we have
\begin{align}\label{rep1}
\tilde{\mu} \| e_{0}^{1/2} \nabla \bm{v} \|_{\bm{\mathcal{W}}}^{2}  & =  \mathbf{w}^{\top} P_{\mathcal{A},1} \mathbf{w}, \quad 
\tilde{\lambda}^{-1} \|e_{0}^{1/2} q_{1} \|_{\mathcal{W}}^{2}  = \bm{q}_{1}^{\top}P_{\mathcal{A},2}  \bm{q}_{1}, \quad \tilde{s}_{0} \|e_{0}^{1/2} q_{2} \|_{\mathcal{W}}^{2}  = \bm{q}_{2}^{\top} P_{\mathcal{A},3} \bm{q}_{2}, 
\end{align}
where $\bm{w} \in \mathbb{R}^{2n_{u}n_{y}}$, $\bm{q}_{1} \in \mathbb{R}^{n_{p}n_{y}}$ and $\bm{q}_{2} \in \mathbb{R}^{n_{p}n_{y}}$ are the associated vectors of coefficients. That is, each of the diagonal blocks of $P_{\mathcal{A}}$ provides a discrete representation of one of the terms in $||| \cdot |||^{2}$ in \eqref{normden}. Applying the action of $P_{\mathcal{A}}^{-1}$ requires only multiple \emph{decoupled} applications of the inverses of the sparse matrices $\mathbb{A}$  and  $\tilde{C}_{0}$. Since $\mathbb{A}$ is a weighted stiffness matrix (a Laplacian matrix if $e_{0}=1$) and $\tilde{C}_{0}$ is a weighted mass matrix, there are many strategies for approximating the actions of $\mathbb{A}^{-1}$  and  $\tilde{C}_{0}^{-1}$ efficiently.

 \subsection{\rbl{Approximation of $\mathcal{S}$}}
Having chosen an approximation $P_{\mathcal{A}}$ to ${\mathcal{A}}$, it would seem quite natural to approximate the Schur complement by $\mathcal{S}_{\textrm{approx}}:= \mathcal{B} P_{\mathcal{A}}^{-1}\mathcal{B}^{\top}+P_{\mathcal{C}}$ where
 \begin{align}
 P_{\mathcal{C}} := & \left(\begin{array}{cc}
I \otimes D_0  & \mathbf{0}\\ \smallskip
\mathbf{0} & \mathbf{0}
 \end{array} \right).
 \end{align} 
However,  this leads to a matrix $ \mathcal{S}_{\textrm{approx}} $ with \emph{dense} blocks.  Instead, we consider the block-diagonal matrix
  \begin{align}\label{PS-def}
 P_{\mathcal{S}}:= \left(\begin{array}{cc} 
 (\alpha^2\tilde{\lambda}^{-1}+\tilde{s}_0) I \otimes \bar{C}_{b}+ I\otimes D_0 & \mathbf{0}\\
\mathbf{0} & (\tilde{\mu}^{-1}+\tilde{\lambda}^{-1})\, I \otimes \bar{C}
\end{array}\right)  =: \left(\begin{array}{cc} 
 P_{\mathcal{S},1} & \bm{0} \\ 
 \bm{0} & P_{\mathcal{S},2}
\end{array}\right)
 \end{align}
where $\bar{C}_{b}$ is the $n_{0} \times n_{0}$ weighted mass matrix associated with the finite element space $\widetilde{V}_{0,p}^{h}$ defined by
\begin{align}
[\bar{C}_{b}]_{r,s}= \int_D e_0(\bm{x})^{-1} \varphi_r(\bm{x}) \varphi_s(\bm{x}) d \bm{x},\quad r,s =1,\ldots,n_{0},
\end{align}
and $\bar{C}$ is the analogous $n_{p} \times n_{p}$ matrix associated with $W^{h}$. For any $q_{F} \in \widetilde{V}_{h,\Lambda}^{0}$ and $q_{T} \in W_{h,\Lambda},$ we have
\begin{align}\label{rep2}
(\alpha^{2} \tilde{\lambda}^{-1}+\tilde{s}_0 ) \| e_{0}^{-1/2} q_{F} \|_{\mathcal{W}}^{2} +  \| \kappa_{0}^{1/2} \nabla q_{F} \|_{\mathcal{W}}^{2}  & =  \mathbf{q}_{F}^{\top} P_{\mathcal{S},1} \mathbf{q}_{F},  \quad
 (\tilde{\mu}^{-1}+\tilde{\lambda}^{-1})  \|e_{0}^{-1/2} q_{T} \|_{\mathcal{W}}^{2}  = \bm{q}_{T}^{\top}P_{\mathcal{S},2}  \bm{q}_{T}, 
\end{align}
where $\bm{q}_{F} \in \mathbb{R}^{n_{0}n_{y}}$ and $\bm{q}_{T}\in \mathbb{R}^{n_{p}n_{y}}$ are the associated vectors of coefficients.  Again, the action of $P_{\mathcal{S}}^{-1}$ can be applied via \emph{decoupled} applications of the inverses of two sparse finite element matrices, namely $(\alpha^2\tilde{\lambda}^{-1} +\tilde{s}_0)\bar{C}_{b}+ D_0$ and the weighted mass matrix $\bar{C}$. 

Combining \eqref{rep1} and \eqref{rep2},  we see that \eqref{disc_norm_rep} holds. As outlined above, bounds for the norm of the underlying preconditioned PDE operator can be derived using the arguments in \cite{mardal2011preconditioning}.  Alternatively, working in the matrix setting, explicit bounds for the eigenvalues of the preconditioned saddlepoint matrix can be derived by following the arguments in Section 4 of \cite{KPSpre}. Using either approach, one obtains bounds that depend on upper and lower bounds for $e_{0}^{-1}E$ and $\kappa_{0}^{-1}\kappa$ (over $D \times \Gamma$), the Korn constant $C_{K}$ and the two discrete inf-sup constants $\gamma, \gamma_{L}$ but \emph{not} on the Poisson ratio $\nu$, the Biot--Willis constant $\alpha$, the rescaled storage coefficient $\tilde{s}_{0}$ or any of the SG-MFEM discretization parameters. Formal proofs are omitted but we now demonstrate the robustness of the preconditioner numerically on test problems.


\section{Numerical results}\label{Numer}
 We first consider two artificial test problems and demonstrate the robustness of the preconditioner $P$  defined by \eqref{Prec_diag}---with $P_{\mathcal{A}}$ and $P_{\mathcal{S}}$ defined as in  \eqref{meanpre} and  \eqref{PS-def}---with respect to the physical and discretization parameters.  We use $\bm{Q}_2-Q_1-Q_1-Q_1-Q_1$ mixed finite elements and choose $S_{\Lambda}$ to be the set of polynomials of total degree $p$ or less in the variables $y_1,\ldots,y_{M_1},z_1,\ldots,z_{M_2}$ on $\Gamma=[-1,1]^{M}$ where $M=M_{1}+M_{2}$. In the third example, we consider a benchmark `footing' problem, similar to examples considered in \cite{botti2019numerical} (a \emph{two-field} stochastic Biot model)  and \cite{oyarzua2016locking} (a three-field \emph{deterministic} Biot model). We apply the proposed SG-MFEM approximation scheme to the new five-field parametric model \eqref{sgos2a11}--\eqref{sgos2a11B} and compute estimates of the mean and variance of the displacement $\bm{u}$ and fluid pressure $p_{F}$ for the nearly incompressible case. All experiments were performed on a MacBook Pro laptop with a modest 16GB of memory and a 2.3GHz Intel Core i5 processor using MATLAB 2019b.


\subsection{Example 1} Here, we choose the spatial domain $D=(0,1)^2$ with $\partial D_{\bm{u}}= [0,1)\times\{0\}\cup \{0\}\times[0,1)$ and $\partial D_{p}= (0,1]\times\{1\}\cup \{1\}\times(0,1]$. The forcing function is chosen to be $\bm{f}=(1,1)^\top$ and $g=0$.  We start by modelling the Young modulus $E$ and  permeability $\kappa$ as spatially uniform. Specifically, 
\begin{align} \label{Ek_ex1}
E(y)= e_0+e_1 y, \quad \kappa(z)=\kappa_0+\kappa_1 z,
\end{align}
where $e_{0},e_{1}, \kappa_{0}, \kappa_{1}$ are constants and $y$ and $z$ are images of independent $U(-1,1)$ random variables. Note that this is a low-dimensional problem with $M_{1}=1=M_{2}$ and $M=2$.  Preconditioned MINRES iteration counts\footnote{The stopping tolerance on the relative residual error in the preconditioned 2-norm was set to $10^{-6}$.} are recorded in Table \ref{tabeff100} for different values of the physical parameters. For the storage coefficient, we follow \cite{lee2017parameter} and set $s_{0}=\alpha^{2}/\lambda$. In this case, the rescaled coefficient $\tilde{s}_{0}$ depends only on $\nu$ and $\alpha$.  We vary $\nu$ and $\alpha$, as well $e_0$ and $\kappa_0$ (the mean values of $E$ and $\kappa$, respectively) and fix $e_{1}=0.1 \times e_{0}$ and $\kappa_{1}=0.1 \times \kappa_{0}$. We also show results for two finite element meshes associated with `level' numbers $\ell=5$ (giving $n_u=1024,n_p=289,n_0=256$) and $l=6$ (giving $n_u=4096,n_p=1089,n_0=1024$). For the parametric approximation, we chose $p=3$ (giving $n_{y}= 10$). Results obtained with $p=4$  (giving $n_{y}= 15$) are also shown in Table \ref{tabeff100_diffp}. We observe that the iteration counts remain bounded when the physical and discretization parameters are varied. This confirms that the proposed preconditioner is robust. 



\begin{table}[ht!]
\caption{MINRES iteration counts and timing in seconds (in parentheses) for varying $\nu$, $\alpha$, $\kappa_0$, $\kappa_1=0.1\times \kappa_0$ and FEM grid level $\ell$ with $p=3$ fixed. In columns 4--6, $e_0=10^5$ and $e_1=10^4$ and in columns 7--9, $e_0=1$ and $e_1=0.1$. }
\label{tabeff100}
\begin{center}
\begin{tabular}{| c | c | c | c | c | c || c |c |c| }
\hline
\multirow{7}{*}{$\kappa_0=1$}&\multirow{4}{*}{$\alpha=1$}&\multicolumn{1}{|c|} {level}&\multicolumn{1}{|c|} {$\nu=.4$} &\multicolumn{1}{|c|}{$\nu=.499$}&\multicolumn{1}{|c||}{$\nu=.49999$} & \multicolumn{1}{|c|}  {$\nu=.4$} &\multicolumn{1}{|c|}{$\nu=.499$}&\multicolumn{1}{|c|}{$\nu=.49999$} \\
\cline{3-9}
& &$l=5$ &$56(0.45)$ &$71(0.58)$ &$71(0.56)$ & $60(0.48)$ &$72(0.58)$ &$71(0.58)$ \\
& &$l=6$ &$56(1.96)$ &$71(2.48)$ &$71(2.46)$ &$60(2.22)$ &$73(2.61)$ &$72(2.55)$ \\
\cline{2-9}
&\multirow{2}{*}{$\alpha=10^{-2}$}&$l=5$ &$55(0.44)$ &$70(0.58)$ &$70(0.58)$ &$56(0.46)$ &$71(0.58)$ &$71(0.59)$ \\
& &$l=6$ &$55(1.94)$ &$71(2.47)$ &$70(2.47)$  &$58(2.02)$ &$71(2.51)$ &$71(2.52)$ \\
\cline{2-9}
&\multirow{2}{*}{$\alpha=10^{-4}$}&$l=5$ &$55(0.45)$ &$70(0.58)$ &$70(0.57)$ &$56(0.46)$ &$70(0.58)$ &$70(0.59)$ \\
& &$l=6$ &$55(1.94)$ &$70(2.47)$ &$70(2.50)$ &$56(1.98)$ &$71(2.5)$ &$71(2.52)$  \\
\cline{1-9}
\multirow{6}{*}{$\kappa_0=10^{-5}$}&\multirow{2}{*}{$\alpha=1$}&$l=5$ &$60(0.48)$ &$72(0.59)$ &$71(0.57)$ &$71(0.63)$ &$70(0.60)$ &$72(0.60)$ \\
& &$l=6$ &$60(2.11)$ &$73(2.56)$ &$72(2.53)$ &$72(2.62)$ &$73(2.65)$ &$73(2.56)$  \\
\cline{2-9}
&\multirow{2}{*}{$\alpha=10^{-2}$}&$l=5$ &$56(0.46)$ &$71(0.59)$ &$71(0.57)$ &$64(0.52)$ &$72(0.59)$ &$71(0.59)$  \\
& &$l=6$ &$58(2.02)$ &$71(2.47)$ &$71(2.51)$ &$64(2.27)$ &$73(2.57)$ &$73(2.61)$ \\
\cline{2-9}
&\multirow{2}{*}{$\alpha=10^{-4}$}&$l=5$ &$56(0.45)$ &$70(0.58)$ &$70(0.56)$  &$58(0.46)$ &$71(0.59)$ &$71(0.58)$ \\
& &$l=6$ &$56(1.94)$ &$71(2.5)$ &$71(2.50)$ &$58(2.05)$ &$71(2.52)$ &$71(2.48)$ \\
\cline{1-9}
\multirow{6}{*}{$\kappa_0=10^{-10}$}&\multirow{2}{*}{$\alpha=1$}&$l=5$ &$71(0.61)$ &$70(0.59)$ &$72(0.63)$ &$71(0.61)$ &$70(0.59)$ &$68(0.60)$ \\
& &$l=6$ &$72(2.61)$ &$73(2.64)$ &$73(2.63)$ &$72(2.63)$ &$70(2.53)$ &$70(2.57)$   \\
\cline{2-9}
&\multirow{2}{*}{$\alpha=10^{-2}$}&$l=5$ &$64(0.52)$ &$72(0.59)$ &$71(0.57)$  &$71(0.64)$ &$70(0.61)$ &$71(0.60)$ \\
& &$l=6$ &$64(2.25)$ &$73(2.57)$ &$73(2.62)$ &$72(2.70)$ &$70(2.56)$ &$73(2.76)$  \\
\cline{2-9}
&\multirow{2}{*}{$\alpha=10^{-4}$}&$l=5$ &$58(0.46)$ &$71(0.58)$ &$71(0.57)$ &$70(0.60)$ &$72(0.58)$ &$71(0.58)$ \\
& &$l=6$ &$58(2.04)$ &$71(2.53)$ &$71(2.51)$  &$70(2.54)$ &$73(2.57)$ &$73(2.55)$ \\
\hline
\end{tabular}
\end{center}
\end{table} 

\begin{table}[ht!]
\caption{MINRES iteration counts and timing in seconds (in parentheses) for varying $\nu$, $\alpha$, $\kappa_0$, $\kappa_1=0.1\times \kappa_0$ and FEM grid level $\ell$  with $p=4$ fixed. In columns 4--6, $e_0=10^5$ and $e_1=10^4$ and in columns 7--9, $e_0=1$ and $e_1=0.1$. }
\label{tabeff100_diffp}
\begin{center}
\begin{tabular}{| c | c | c | c | c | c || c |c |c| }
\hline
\multirow{7}{*}{$\kappa_0=10^{-10}$}&\multirow{4}{*}{$\alpha=1$}&\multicolumn{1}{|c|} {level}&\multicolumn{1}{|c|} {$\nu=.4$} &\multicolumn{1}{|c|}{$\nu=.499$}&\multicolumn{1}{|c||}{$\nu=.49999$} & \multicolumn{1}{|c|}  {$\nu=.4$} &\multicolumn{1}{|c|}{$\nu=.499$}&\multicolumn{1}{|c|}{$\nu=.49999$} \\
\cline{2-9}
& &$l=5$ &$72 ( 0.78)$ &$70( 0.74)$ &$73(0.77 )$ &$72 (0.77 )$ &$70 (0.74 )$ &$70 (0.74 )$ \\
& &$l=6$ &$72 (3.75 )$ &$ 73(3.83 )$ &$73 (3.81)$ &$ 72(3.72 )$ &$70 ( 3.63)$ &$70 (3.66 )$   \\
\cline{2-9}
&\multirow{2}{*}{$\alpha=10^{-2}$}&$l=5$ &$64(0.70 )$ &$73 (0.78 )$ &$ 73(0.83 )$  &$72 (0.76 )$ &$70 (0.73 )$ &$71 (0.73 )$ \\
& &$l=6$ &$ 64(3..37 )$ & $73 (3.82 )$ & $ 73(3.79 )$ &$72 (3.72 )$ & $ 70(3.68 )$ &$  73(3.79 )$  \\
\cline{2-9}
&\multirow{2}{*}{$\alpha=10^{-4}$}&$l=5$ &$ 58(0.60 )$ &$ 71(0.76 )$ &$ 71( 0.76)$ &$ 70(0.75 )$ &$ 73(0.77 )$ &$ 73(0.77 )$ \\
& &$l=6$ &$58 (3.08 )$ &$ 73( 3.77)$ &$ 72(3.74 )$ &$ 71(3.83 )$ &$74 (3.85 )$ & $73 (3.79 )$ \\
\hline
\end{tabular}
\end{center}
\end{table} 

\subsection{Example 2} Next, we consider $D=(-1,1)^2$ with $\partial D_{\bm{u}}= [-1,1)\times\{-1\}\cup \{-1\}\times[-1,1)$ and $\partial D_{p}= (-1,1]\times\{1\}\cup \{1\}\times(-1,1]$. Again, we choose $\bm{f}=(1,1)^\top$ and $g=0$. This time, the Young modulus $E$ and hydraulic conductivity $\kappa$ are modelled as
\begin{align} \label{uncertainE}
E(\bm{x},\bm{y})= e_0 +  \sigma_E  \sum_{m=1}^{M_{1}}\sqrt{\lambda_m}\varphi_m(\bm{x})y_m, \qquad
\kappa(\bm{x},\bm{z})= \kappa_0 +  \sigma_{\kappa} \sum_{m=1}^{M_2}\sqrt{\lambda_m}\varphi_m(\bm{x})z_m,
\end{align}
where $y_{m}, z_{m}$ are images of independent $U(-\sqrt{3},\sqrt{3})$ random variables, $\sigma_E$ and $\sigma_{\kappa}$ are the standard deviations of the respective random fields, $e_{0}$ and $\kappa_{0}$ are the (constant) means, and $\{(\lambda_m,\varphi_m)\}$ are the eigenpairs 
of the integral operator associated with the covariance kernel
\begin{align}
C(\bm{x},\bm{x}')= \exp\left(-\frac{1}{2} ||\bm{x}-\bm{x}' ||_{1} \right), \quad\bm{x},\bm{x}'\in D.
\end{align}
Note that this means that $E$ and $\kappa$ have the covariance functions $\sigma^{2}_EC(\bm{x}, \bm{x}')$ and $\sigma^{2}_{\kappa}C(\bm{x}, \bm{x}')$, respectively. We do not have to choose the same spatial covariance kernel for both fields (indeed, this would be unphysical in realistic applications); we do so here for simplicity only to test the performance of the preconditioner.  We define the storage coefficient as described in Example 1 and vary the physical parameters  $\nu$, $\alpha$ and $\kappa_{0}$ as before with $\sigma_{\kappa}=0.1\times \kappa_0$. For the Young modulus, we fix the mean and standard deviation to be $e_0=10^5$ and $\sigma_{E}=0.1\times e_{0}$. In Table \ref{tabeff3+3} we show results for two cases: $M_1=M_2=3$ (giving $n_y=84$) and $M_1=M_2=5$ (giving $n_y=286)$.  Again, we consider two finite element meshes with level numbers $\ell=5,6$  and now fix $p=3$ for the parametric approximation. The MINRES iteration counts remain bounded as the discretization and physical parameters are varied, confirming once again that the preconditioner is robust. 

\begin{table}
\caption{MINRES iteration counts and timing in seconds (in parentheses) for varying $\nu$, $\alpha$, $\kappa_0$,  $\sigma_{\kappa}=0.1\times \kappa_0$ and FEM grid level $\ell$, with $e_0=10^5$, $\sigma_{E}=0.1\times e_{0}$ and $p=3$ fixed. In columns 4--6, $M_{1}=3, M_{2}=3$ and in columns 7--9, $M_{1}=5, M_{2}=5$. }
\label{tabeff3+3}
\begin{center}
\begin{tabular}{| c | c | c | c | c | c || c| c| c|  }
\hline
\multirow{8}{*}{{\small $\kappa_0=$}}&\multirow{4}{*}{{\small $\alpha=1$}}&\multicolumn{1}{|c|} {level} &\multicolumn{1}{|c|} {\small $\nu=.4$}&\multicolumn{1}{|c|}{\small $\nu=.499$}&\multicolumn{1}{|c||}{{\small $\nu=.49999$}} &\multicolumn{1}{|c|} {\small $\nu=.4$}&\multicolumn{1}{|c|}{\small $\nu=.499$}&\multicolumn{1}{|c|}{{\small $\nu=.49999$}} \\
\cline{3-9}
& &$l=5$ &$79(4.23)$ &$98(5.33)$ &$96(5.31)$ &$81(16.0)$ &$99(18.9)$ &$99(18.5)$ \\
& &$l=6$ &$81(22.7)$ &$98(27.3)$ &$98(27.5)$ &$82(90.2)$ &$101(113.0)$ &$101(113.8)$ \\
\cline{2-9}
&\multirow{2}{*}{{$\alpha=10^{-2}$}}&$l=5$ &$79(4.33)$ &$96(5.29)$ &$95(5.18)$ &$80(14.6)$ &$99(18.4)$ &$97(17.9)$ \\
& &$l=6$ &$79(23.1)$ &$98(27.4)$ &$97(27.1)$ &$80(86.9)$ &$99(110.6)$ &$99(107.5)$  \\
\cline{2-9}
1 &\multirow{2}{*}{{$\alpha=10^{-4}$}}&$l=5$ &$77(4.12)$ &$95(5.15)$ &$95(5.12)$ &$78(15.1)$ &$97(17.7)$ &$96(17.4)$ \\
& &$l=6$ &$78(22.1)$ &$97(27.1)$ &$96(27.0)$  &$79(87.3)$ &$99(108.3)$ &$98(106.7)$ \\
\cline{1-9}
\multirow{6}{*}{{\small $\kappa_0=$}}&\multirow{2}{*}{{\small $\alpha=1$}}&$l=5$ &$85(4.88)$ &$99(5.42)$ &$98(5.35)$ &$86(15.7)$ &$100(18.3)$ &$100(18.3)$  \\
& &$l=6$ &$85(24.3)$ &$100(29.0)$ &$100(29.1)$ &$87(94.0)$ &$102(112.3)$ &$101(109.0)$  \\
\cline{2-9}
&\multirow{2}{*}{{$\alpha=10^{-2}$}}&$l=5$ &$80(4.38)$ &$98(5.28)$ &$96(5.37)$  &$81(14.9)$ &$100(18.8)$ &$99(18.8)$ \\
 & &$l=6$ &$81(22.7)$ &$100(27.9)$ &$99(27.7)$ &$83(89.6)$ &$101(109.8)$ &$101(109.5)$ \\
\cline{2-9}
&\multirow{2}{*}{$\alpha=10^{-4}$}&$l=5$ &$79(4.34)$ &$96(5.4)$ &$95(5.11)$ &$80(15.2)$ &$99(18.1)$ &$97(18.0)$  \\
$ 10^{-5}$ & &$l=6$ &$79(22.5)$ &$98(27.3)$ &$97(27.3)$ &$81(87.4)$ &$99(107.1)$ &$99(107.1)$ \\
\cline{1-9}
\multirow{6}{*}{{\small $\kappa_0= $}}&\multirow{2}{*}{{\small $\alpha=1$}}&$l=5$ &$96(5.58)$ &$91(5.22)$ &$98(5.71)$ &$97(19.7)$ &$93(17.9)$ &$100(20.6)$ \\
 & &$l=6$ &$97(28.7)$ &$96(27.3)$ &$100(29.3)$ &$98(114.4)$ &$97(113.8)$ &$102(119.3)$ \\
\cline{2-9}
&\multirow{2}{*}{{$\alpha=10^{-2}$}}&$l=5$ &$93(5.20)$ &$99(5.34)$ &$98(5.32)$ &$95(17.4)$ &$100(18.6)$ &$100(18.3)$ \\
  & &$l=6$ &$94(26.9)$ &$100(28.0)$ &$100(28.2)$ &$95(103.4)$ &$102(111.5)$ &$102(111.6)$  \\
\cline{2-9}
&\multirow{2}{*}{{$\alpha=10^{-4}$}}&$l=5$ &$81(4.40)$ &$98(5.28)$ &$96(5.19)$  &$81(14.8)$ &$100(18.4)$ &$99(18.0)$ \\
$10^{-10}$ & &$l=6$ &$81(22.8)$ &$100(28.1)$ &$99(27.7)$ &$83(90.3)$ &$101(109.6)$ &$101(109.6)$ \\
\hline
\end{tabular}
\end{center}
\end{table}




\subsection{Example 3}
\begin{figure}
\includegraphics[width=2.7in]{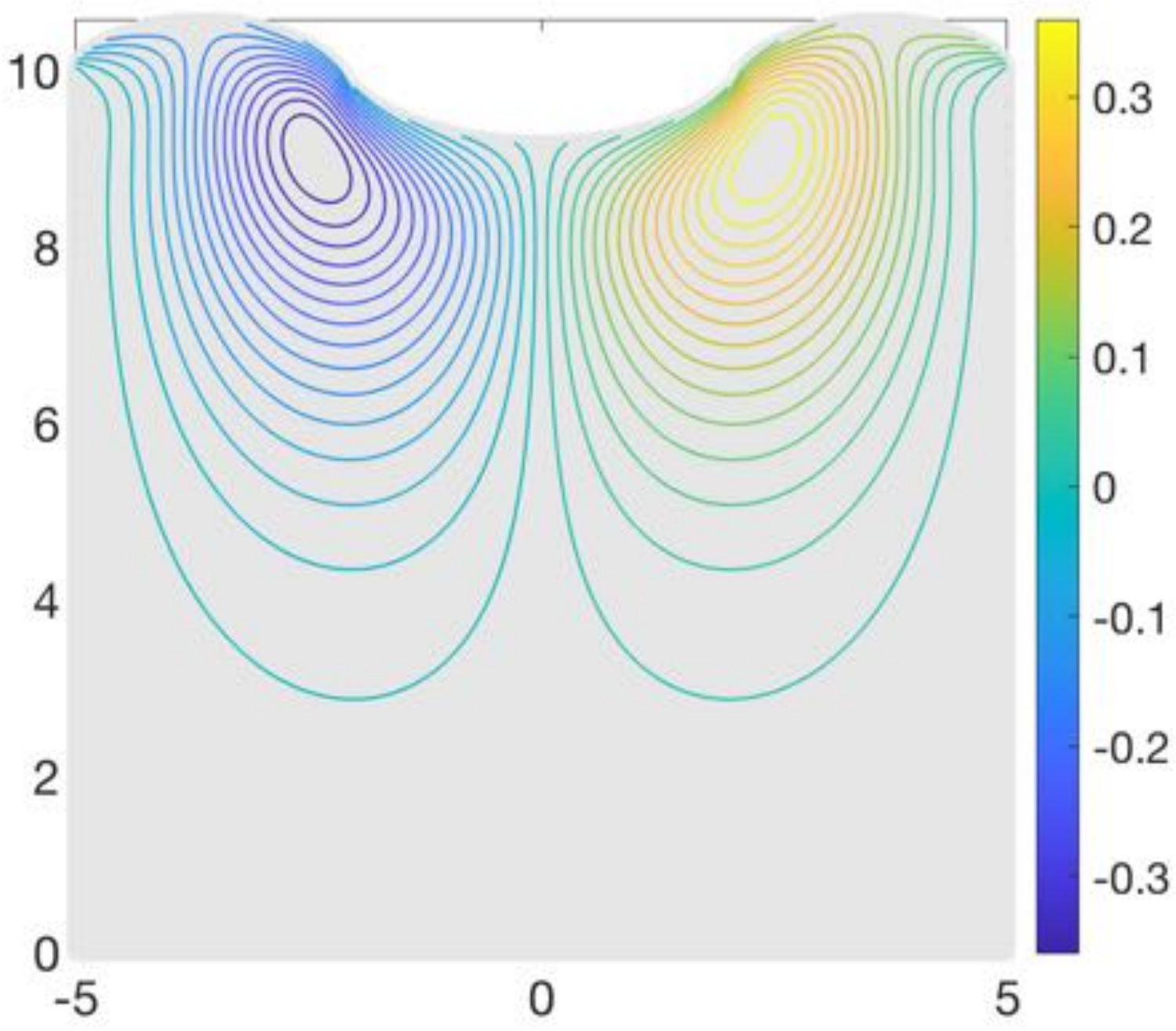}
\includegraphics[width=2.7in]{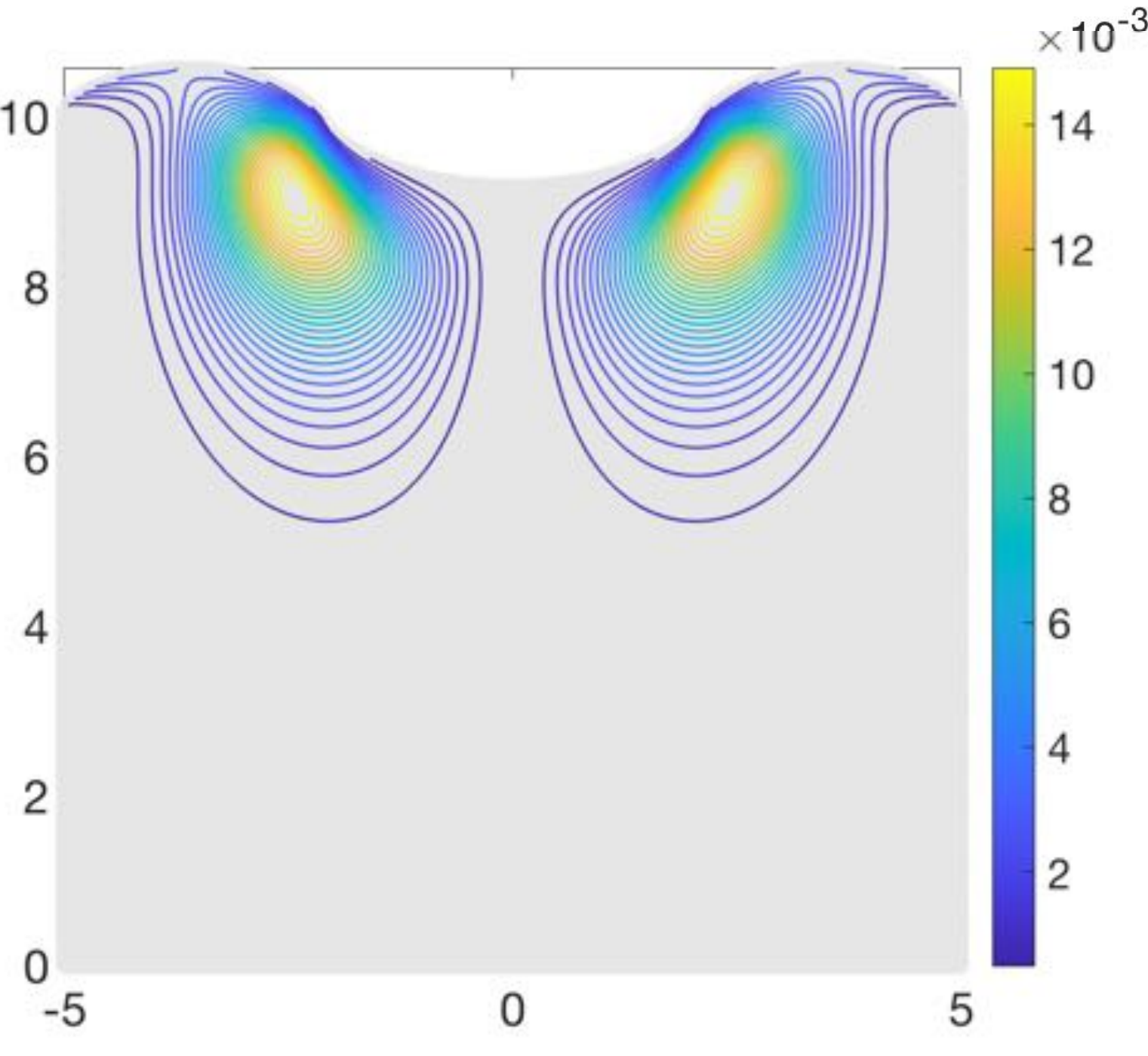}
\includegraphics[width=2.7in]{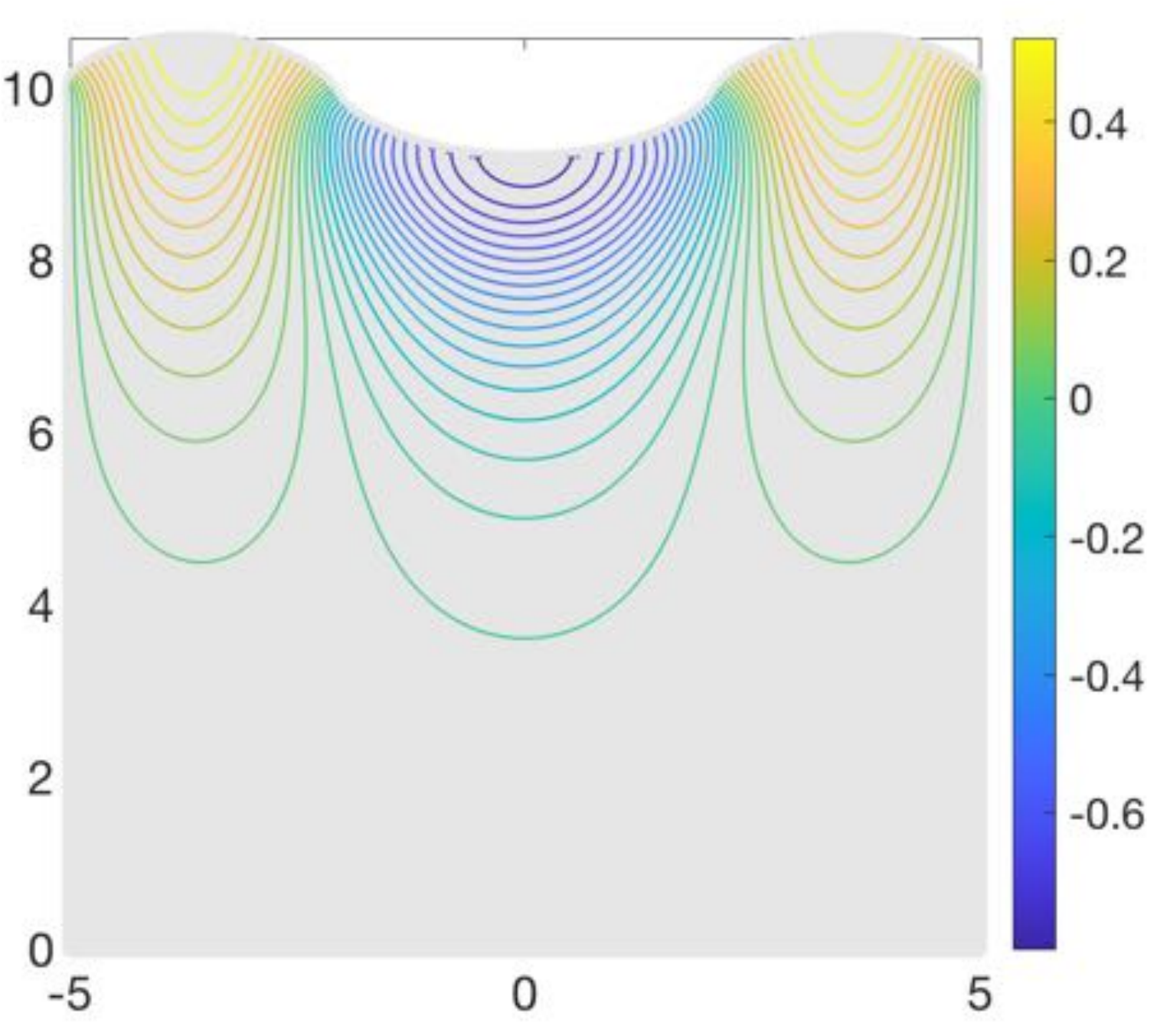}
\includegraphics[width=2.7in]{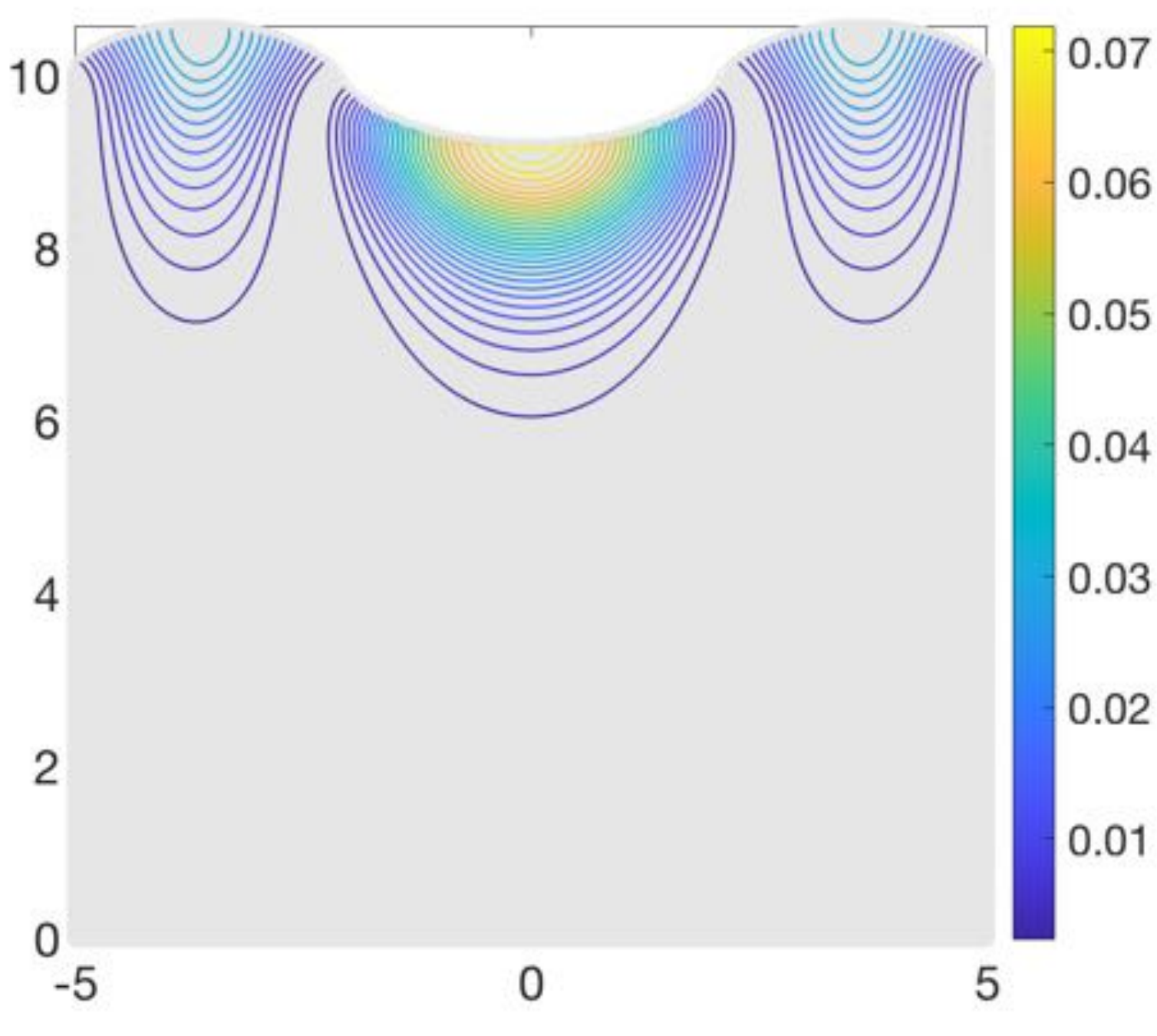}
\includegraphics[width=2.7in]{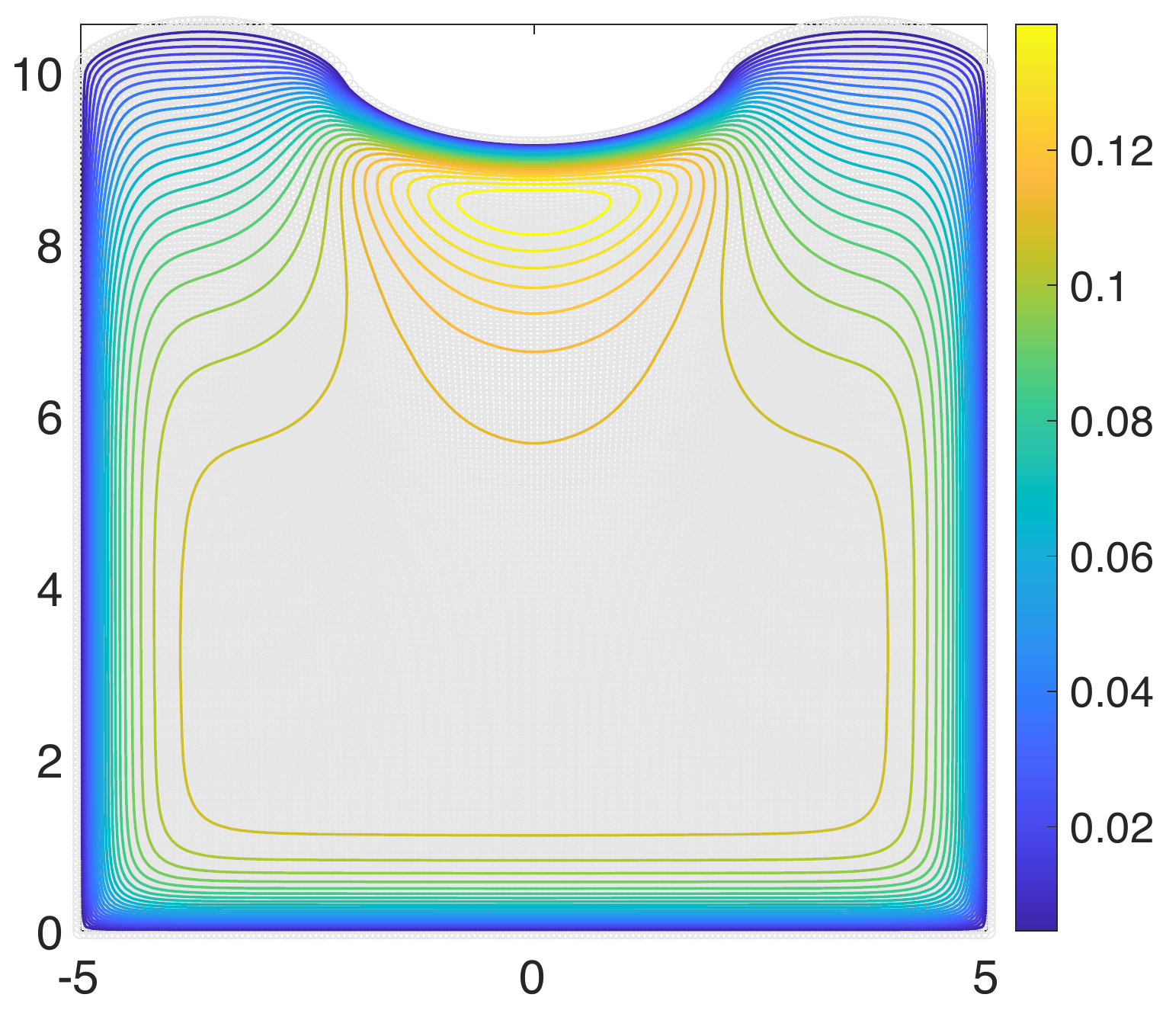}
\hspace{0.9in}
\includegraphics[width=2.7in]{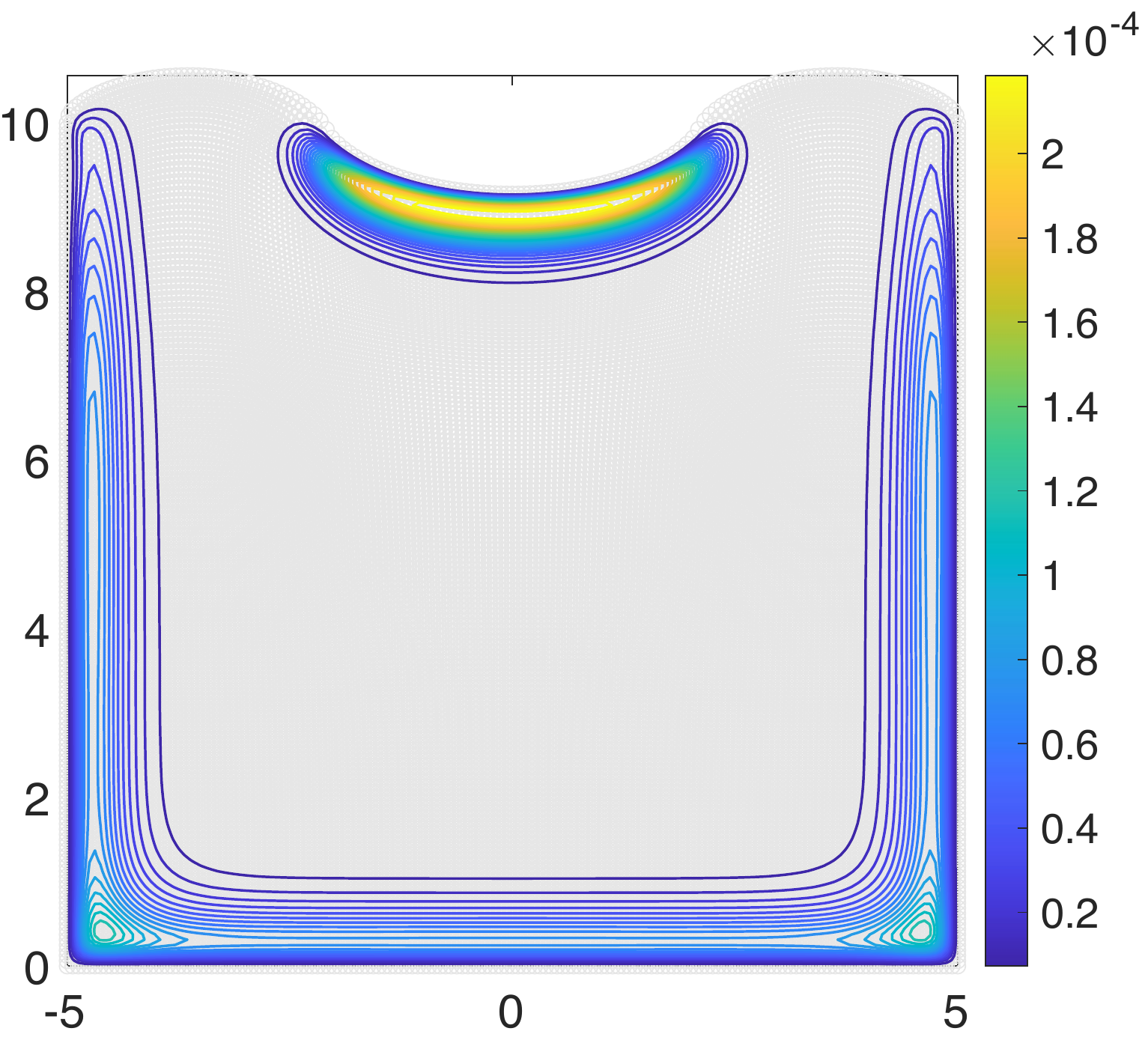}
\caption{\label{Biot_pic} Mean (left) and variance (right) of the SG-MFEM approximation to the horizontal displacement $u_{1}$ (top), the vertical displacement $u_{2}$ (middle) and the fluid pressure $p_{F}$ (bottom) computed for Example 3 using $\bm{Q}_2-Q_1-Q_1-Q_1-Q_1$ mixed finite elements with mesh level $\ell=5$ and polynomials of total degree $p=4$ or less. The physical parameters are $\nu=0.4995$, $\alpha=0.1$ and $\tilde{s}_0=30.$}
\end{figure}

\begin{figure}
\includegraphics[width=2.7in]{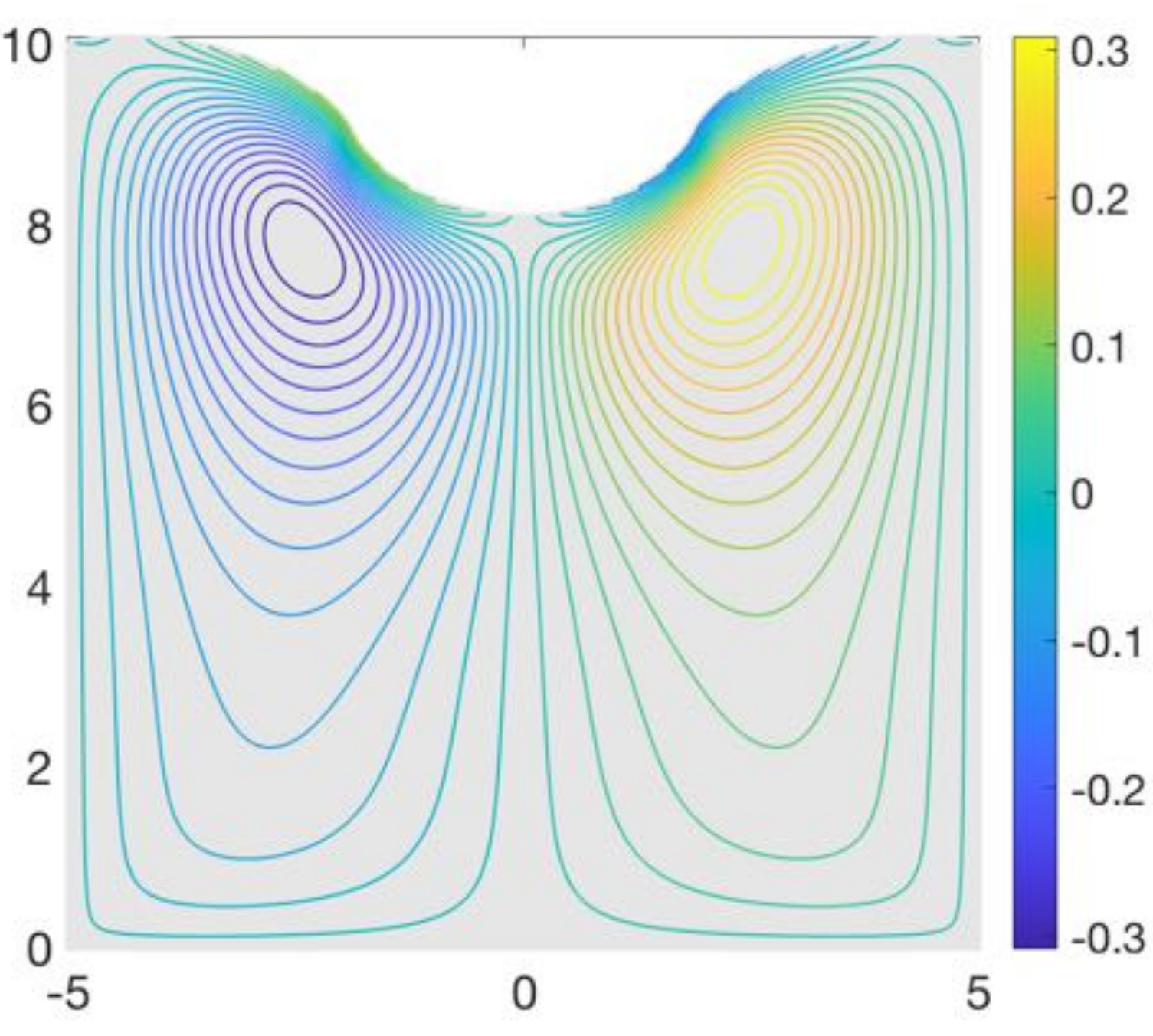}
\includegraphics[width=2.7in]{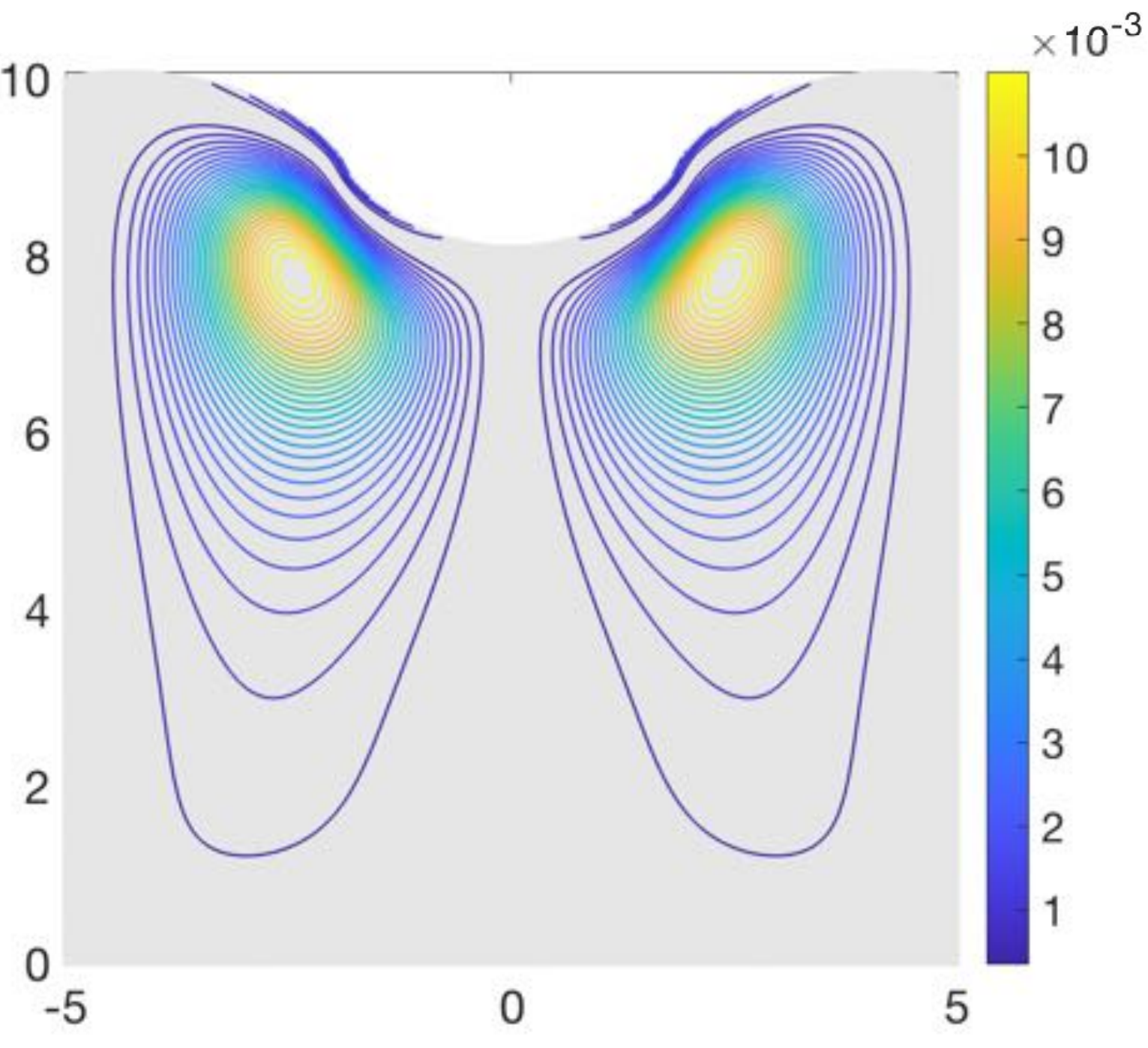}
\includegraphics[width=2.7in]{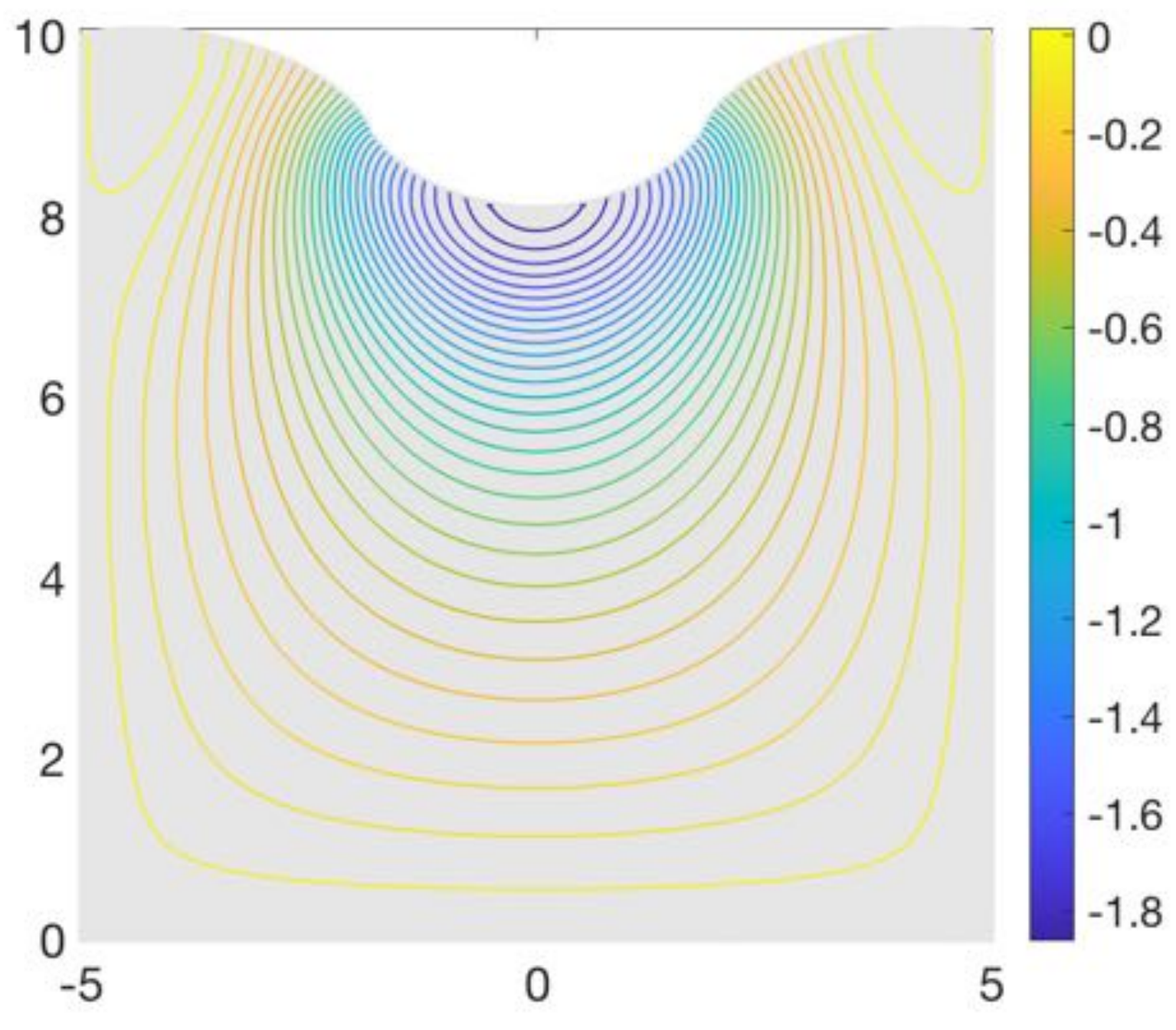}
\includegraphics[width=2.7in]{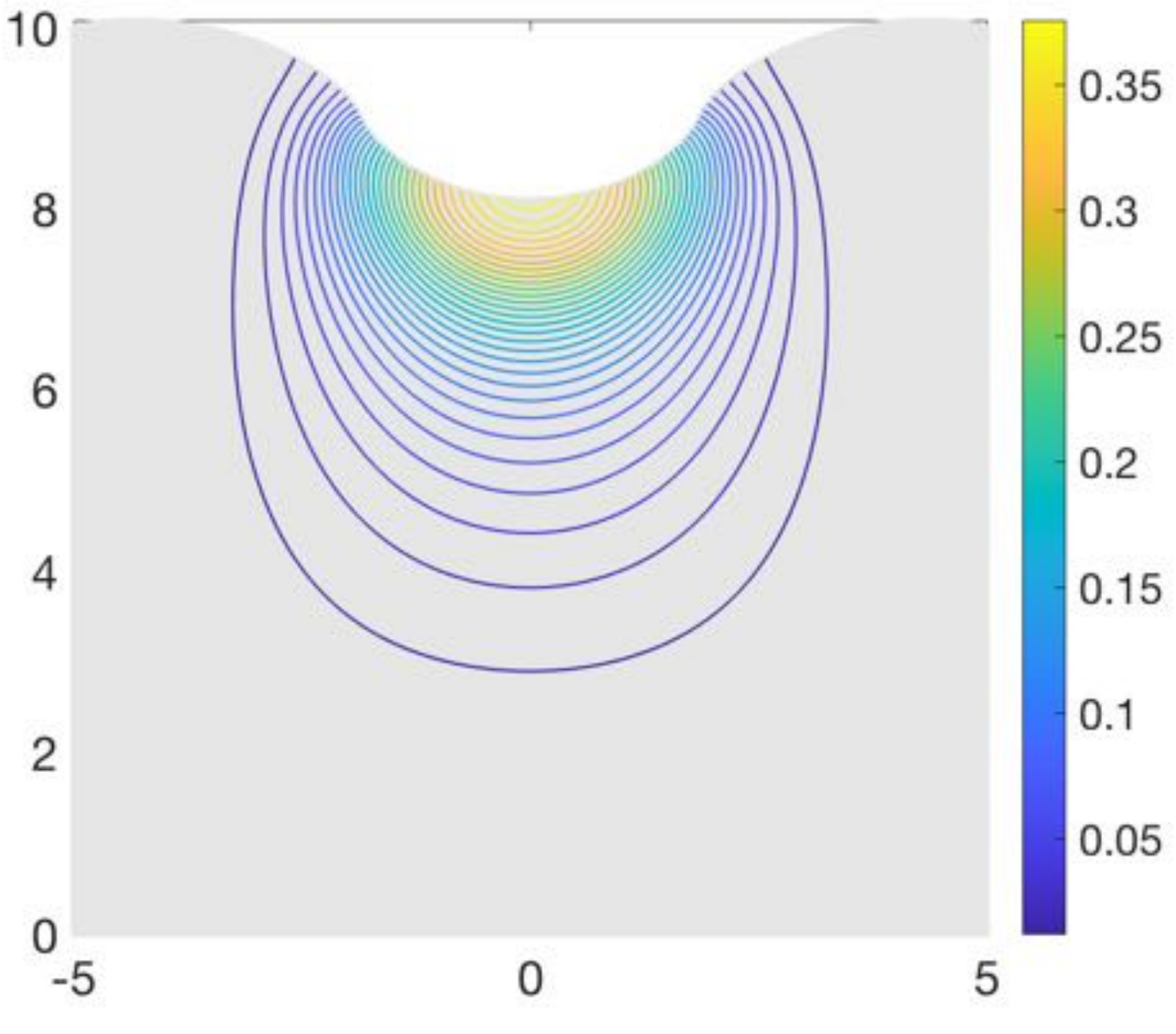}
\includegraphics[width=2.7in]{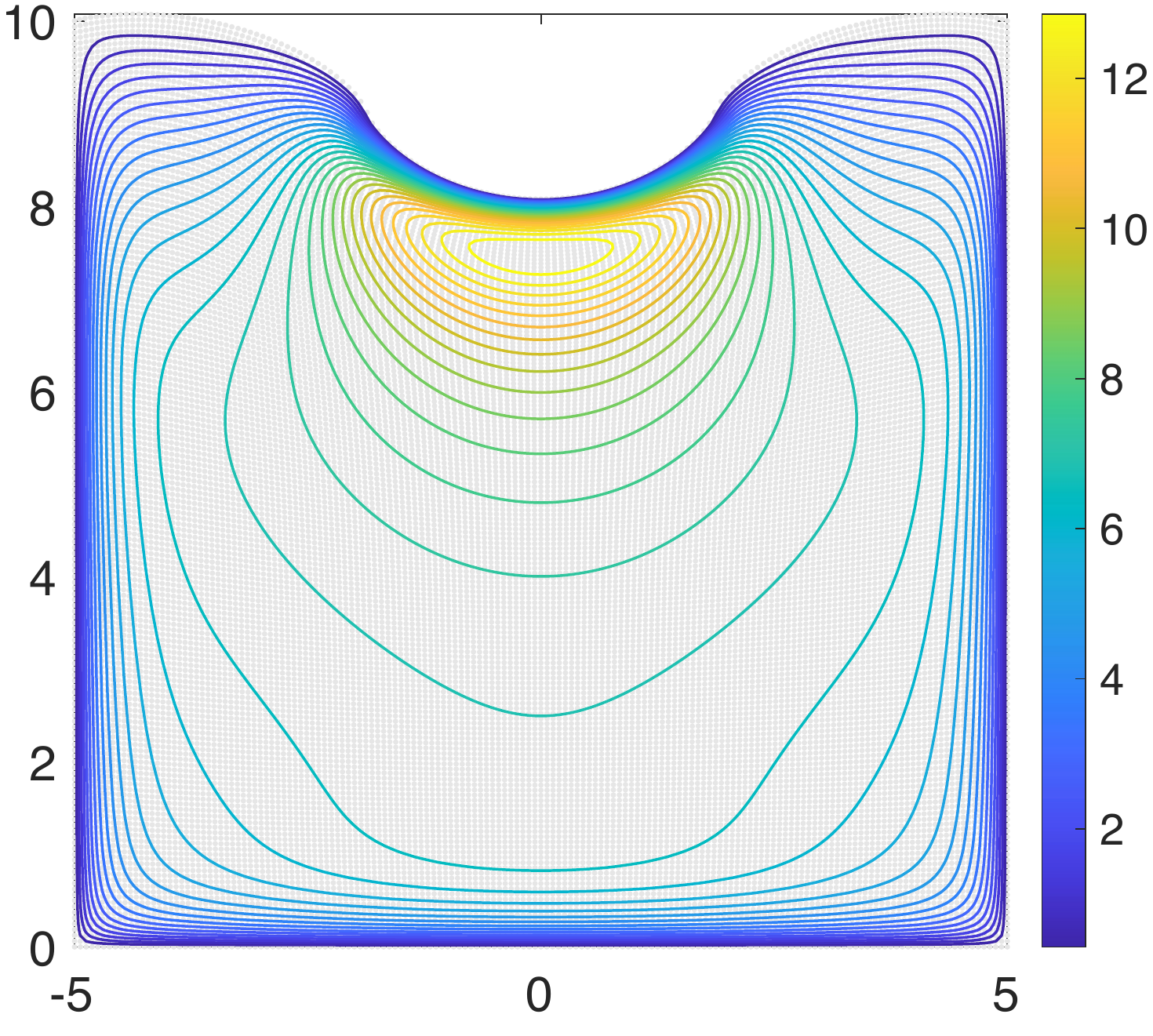}
\hspace{0.9in}
\includegraphics[width=2.7in]{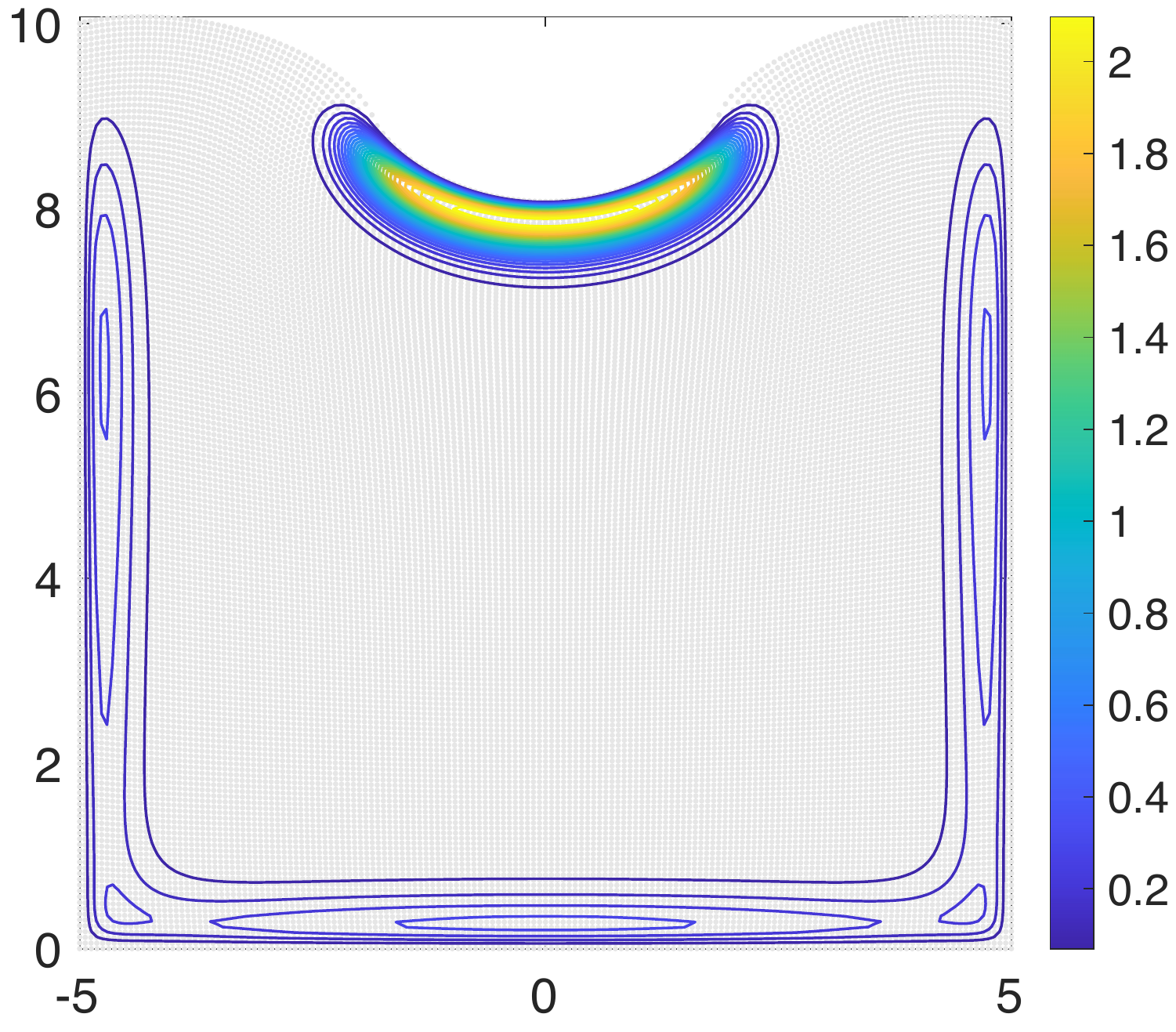}
\caption{\label{Biot_pic2} Mean (left) and variance (right) of the SG-MFEM approximation to the horizontal displacement $u_{1}$ (top), the vertical displacement $u_{2}$ (middle) and the fluid pressure $p_{F}$ (bottom) computed for Example 3 using $\bm{Q}_2-Q_1-Q_1-Q_1-Q_1$ mixed finite elements with mesh level $\ell=5$ and polynomials of total degree $p=4$ or less. The physical parameters are $\nu=0.45$, $\alpha=0.1$ and $\tilde{s}_0=30.$}
\end{figure}

\textbf{(Footing Problem)} We consider a block of porous soil saturated with fluid that has a load of intensity $h$ applied in the downward vertical direction on a portion of the top side. The spatial domain is chosen to be $D=(-5,5)\times(0,10)$ and we define $\partial D_{p}:=(-5,5)\times \{10\}$ and $\partial D_{\bm{u}}= \partial D \setminus \partial D_{p}$. We also choose $\bm{f}=(0,0)^\top$ and $g=0$. The Young modulus $E$ and hydraulic conductivity $\kappa$ are modelled as in \eqref{Ek_ex1} in Example 1 and so the parameter domain is $\Gamma:=[-1,1]\times [-1,1]$. For the boundary conditions, we set $p_F= 0$ on $\partial D \times \Gamma$,  $\bm{u}= \bm{0}$ on $\partial D_{\bm{u}} \times \Gamma$ and  $\sigmab\bm{n}=(0,-h)$ on $\partial D_p \times \Gamma$ where
\begin{align}
h=\begin{cases}
1.5 \times 10^{4} & x_{1} \in[-2,2],\\
0 & \textrm{otherwise}.
\end{cases}
\end{align}
For the physical parameters, we follow the deterministic example in \cite{oyarzua2016locking} and choose $\alpha=0.1$, $\tilde{s}_0=30$ and $\nu=0.4995$ (so the material is nearly incompressible) and we choose the mean values of $E$ and $\kappa$ to be $e_0=3\times 10^{4}$ and $\kappa_0=10^{-4}$.  For the stochastic part, we then choose $e_1=0.5\times e_0$ and $\kappa_1=0.5 \times \kappa_0.$ The mean and variance of the SG-MFEM approximation to the components of the displacement and the fluid pressure computed with mesh level $\ell=5$ (here,  giving $n_{x}=307,204$) and $p=4$ (giving $n_{y}=15$) are plotted in Figure \ref{Biot_pic}. The mean solution fields are qualitatively similar to the deterministic results shown in Figure 3 in \cite{oyarzua2016locking}. A similar time-dependent footing problem was also considered in \cite{botti2019numerical}, although the approximation scheme considered there is not suitable in the nearly incompressible case. For comparison, results obtained with $\nu=0.45$ are also plotted in Figure \ref{Biot_pic2}. Observe that even when close to the incompressible limit, there no spurious oscillations.  

 \section{Summary}  We presented a new parameter-robust approximation scheme and a solver for a linear poroelasticity model with uncertain inputs. Specifically, we introduced a five-field Biot model with uncertain Young modulus $E$ and hydraulic conductivity field $\kappa$ and rigorously analyzed the well-posedness and stability of the associated weak problem with respect to a weighted norm.  We explained how to discretize the weak problem using a locking-free stochastic Galerkin mixed finite element method and then constructed a parameter-robust preconditioner for the associated linear systems which is informed by the weighted norm with respect to which the approximation is provably stable. Numerical results demonstrate the robustness of the preconditioner with respect to three key physical parameters, as well as the SG-MFEM discretization parameters.  Crucially, both the approximation scheme and solver are robust when the Poisson ratio $\nu \to 1/2$.  As far as we are aware, this is the first attempt to develop an approximation scheme for stochastic Biot-type problems that is robust in the incompressible limit. 
 
 However, there are clearly many possible avenues for future work.  So far, we have not examined transient problems. We also focused on uncertainty in the spatially varying inputs $E$ and $\kappa$. Other physical parameters could also be treated as uncertain. Although we achieved a preconditioning scheme that is robust with respect to the Poisson ratio, the Biot--Willis constant and the storage coefficient, it is not robust with respect to variations in $E_{\max}/e_{0}^{\min}$ and $\kappa_{\max}/\kappa_{0}^{\min}$. This is because the preconditioner only incorporates the leading terms $e_{0}$ and $\kappa_{0}$ from the representations of $E$ and $\kappa$. This leads to a block-diagonal preconditioner, which can therefore be implemented efficiently. However, convergence will deteriorate if the standard deviations of $E$ and $\kappa$ are large relative to their mean values.  Depending on the computational resources available, more sophisticated solvers may need to be developed for more challenging cases. Finally, now that a discretization scheme and solver are in place, a posteriori error estimation needs to be explored. Galerkin approximation offers a very natural framework for doing this rigorously and recent work \cite{khan2018error3field} for linear elasticity problems provides a suitable starting point.
 
%
 %


\bibliography{kps_upd}

\begin{thebibliography}{10}

\bibitem{berger2015stabilized}
Lorenz Berger, Rafel Bordas, David Kay, and Simon Tavener.
\newblock Stabilized lowest-order finite element approximation for linear
  three-field poroelasticity.
\newblock {\em SIAM Journal on Scientific Computing}, 37(5):A2222--A2245, 2015.

\bibitem{botti2019numerical}
Michele Botti, Daniele~A Di~Pietro, Olivier~Le Ma{\^\i}tre, and Pierre Sochala.
\newblock Numerical approximation of poroelasticity with random coefficients
  using polynomial chaos and hybrid high-order methods.
\newblock {\em arXiv preprint arXiv:1903.11885}, 2019.

\bibitem{chang1985uncertainty}
Ching~S Chang.
\newblock Uncertainty of one-dimensional consolidation analysis.
\newblock {\em Journal of geotechnical engineering}, 111(12):1411--1424, 1985.

\bibitem{Cohen}
Albert Cohen, Ronald DeVore, and Christoph Schwab.
\newblock Analytic regularity and polynomial approximation of parametric and
  stochastic elliptic {PDE}s.
\newblock {\em Analysis and Applications}, 09(01):11--47, 2011.

\bibitem{coussy2004poromechanics}
Olivier Coussy.
\newblock {\em Poromechanics}.
\newblock John Wiley \& Sons, 2004.

\bibitem{Adam_ML}
Adam~J. Crowder, Catherine~E. Powell, and Alex Bespalov.
\newblock Efficient adaptive multilevel stochastic {G}alerkin approximation
  using implicit a posteriori error estimation.
\newblock {\em SIAM Journal on Scientific Computing}, 41(3):A1681--A1705, 2019.

\bibitem{darrag1993consolidation}
AA~Darrag and MA~El~Tawil.
\newblock The consolidation of soils under stochastic initial excess pore
  pressure.
\newblock {\em Applied mathematical modelling}, 17(11):609--612, 1993.

\bibitem{delgado2015stochastic}
Paul Delgado and Vinod Kumar.
\newblock A stochastic {G}alerkin approach to uncertainty quantification in
  poroelastic media.
\newblock {\em Applied Mathematics and Computation}, 266:328--338, 2015.

\bibitem{detournay1993fundamentals}
Emmanuel Detournay and Alexander H-D Cheng.
\newblock Fundamentals of poroelasticity.
\newblock In {\em Analysis and design methods}, pages 113--171. Elsevier, 1993.

\bibitem{EIGEL}
Martin Eigel, Claude~Jeffrey Gittelson, Christoph Schwab, and Elmar Zander.
\newblock Adaptive stochastic {G}alerkin {FEM}.
\newblock {\em Computer Methods in Applied Mechanics and Engineering}, 270:247
  -- 269, 2014.

\bibitem{HDA}
Howard Elman, David Silvester, and Andy Wathen.
\newblock {\em Finite Elements and Fast Iterative Solvers: with Applications in
  Incompressible Fluid Dynamics}.
\newblock Oxford University Press, Oxford, UK, 2014.
\newblock Second Edition, xiv+400 pp. ISBN: 978-0-19-967880-8.

\bibitem{ernst2009efficient}
Oliver~G. Ernst, Catherine~E. Powell, David~J. Silvester, and Elisabeth
  Ullmann.
\newblock Efficient solvers for a linear stochastic {G}alerkin mixed
  formulation of diffusion problems with random data.
\newblock {\em SIAM Journal on Scientific Computing}, 31(2):1424--1447, 2009.

\bibitem{frias2004stochastic}
Diego~G Frias, M{\'a}rcio~A Murad, and Felipe Pereira.
\newblock Stochastic computational modelling of highly heterogeneous
  poroelastic media with long-range correlations.
\newblock {\em International Journal for Numerical and Analytical Methods in
  Geomechanics}, 28(1):1--32, 2004.

\bibitem{relax_prec}
Matteo Frigo, Nicola Castelletto, and Massimiliano. Ferronato.
\newblock A relaxed physical factorization preconditioner for mixed finite
  element coupled poromechanics.
\newblock {\em SIAM Journal on Scientific Computing}, 41(4):B694--B720, 2019.

\bibitem{haga2012causes}
Joachim~Berdal Haga, Harald Osnes, and Hans~Petter Langtangen.
\newblock On the causes of pressure oscillations in low-permeable and
  low-compressible porous media.
\newblock {\em International Journal for Numerical and Analytical Methods in
  Geomechanics}, 36(12):1507--1522, 2012.

\bibitem{khan2018error3field}
Arbaz Khan, Alex Bespalov, Catherine~E. Powell, and David~J. Silvester.
\newblock Robust a posteriori error estimation for stochastic {G}alerkin
  formulations of parameter-dependent linear elasticity equations.
\newblock {\em arXiv preprint}, 2018.
\newblock {\tt https://arxiv.org/abs/1810.07440}.

\bibitem{KPSpre}
Arbaz Khan, Catherine~E Powell, and David~J Silvester.
\newblock Robust preconditioning for stochastic {G}alerkin formulations of
  parameter-dependent nearly incompressible elasticity equations.
\newblock {\em SIAM Journal on Scientific Computing}, 41(1):A402--A421, 2019.

\bibitem{klawonn1998optimal}
Axel Klawonn.
\newblock An optimal preconditioner for a class of saddle point problems with a
  penalty term.
\newblock {\em SIAM Journal on Scientific Computing}, 19(2):540--552, 1998.

\bibitem{lee2017parameter}
Jeonghun~J Lee, Kent-Andre Mardal, and Ragnar Winther.
\newblock Parameter-robust discretization and preconditioning of {B}iot's
  consolidation model.
\newblock {\em SIAM Journal on Scientific Computing}, 39(1):A1--A24, 2017.

\bibitem{lee2019mixed}
Jeonghun~J Lee, Eleonora Piersanti, K-A Mardal, and Marie~E Rognes.
\newblock A mixed finite element method for nearly incompressible
  multiple-network poroelasticity.
\newblock {\em SIAM Journal on Scientific Computing}, 41(2):A722--A747, 2019.

\bibitem{lewis1998finite}
R~W Lewis and Schrefler~B A.
\newblock {\em The Finite Element Method in the Static and Dynamic Deformation
  and Consolidation of Porous Media-RW Lewis and BA Schrefler}.
\newblock Wiley, Chichester, UK, 1998.

\bibitem{mardal2011preconditioning}
Kent-Andre Mardal and Ragnar Winther.
\newblock Preconditioning discretizations of systems of partial differential
  equations.
\newblock {\em Numerical Linear Algebra with Applications}, 18(1):1--40, 2011.

\bibitem{merxhani2016introduction}
Andi Merxhani.
\newblock An introduction to linear poroelasticity.
\newblock {\em arXiv preprint arXiv:1607.04274}, 2016.

\bibitem{murad1992improved}
M{\'a}rcio~A Murad and Abimael~FD Loula.
\newblock Improved accuracy in finite element analysis of {B}iot's
  consolidation problem.
\newblock {\em Computer Methods in Applied Mechanics and Engineering},
  95(3):359--382, 1992.

\bibitem{murad1994stability}
M{\'a}rcio~A Murad and Abimael~FD Loula.
\newblock On stability and convergence of finite element approximations of
  {B}iot's consolidation problem.
\newblock {\em International Journal for Numerical Methods in Engineering},
  37(4):645--667, 1994.

\bibitem{murad1996asymptotic}
M{\'a}rcio~A Murad, Vidar Thom{\'e}e, and Abimael~FD Loula.
\newblock Asymptotic behavior of semidiscrete finite-element approximations of
  {B}iot's consolidation problem.
\newblock {\em SIAM journal on numerical analysis}, 33(3):1065--1083, 1996.

\bibitem{nishimura2002consolidation}
Shin-ichi Nishimura, Kiyoshi Shimada, and Hiroaki Fujii.
\newblock Consolidation inverse analysis considering spatial variability and
  non-linearity of soil parameters.
\newblock {\em Soils and foundations}, 42(3):45--61, 2002.

\bibitem{oyarzua2016locking}
Ricardo Oyarz{\'u}a and Ricardo Ruiz-Baier.
\newblock Locking-free finite element methods for poroelasticity.
\newblock {\em SIAM Journal on Numerical Analysis}, 54(5):2951--2973, 2016.

\bibitem{phillips2009overcoming}
Phillip~Joseph Phillips and Mary~F Wheeler.
\newblock Overcoming the problem of locking in linear elasticity and
  poroelasticity: an heuristic approach.
\newblock {\em Computational Geosciences}, 13(1):5--12, 2009.

\bibitem{powell_elman}
Catherine~E. Powell and Howard~C. Elman.
\newblock Block-diagonal preconditioning for spectral stochastic finite-element
  systems.
\newblock {\em IMA Journal of Numerical Analysis}, 29(2):350--375, 2009.

\bibitem{MultiRB}
Catherine~E. Powell, David Silvester, and Valeria Simoncini.
\newblock An efficient reduced basis solver for stochastic {G}alerkin matrix
  equations.
\newblock {\em SIAM Journal on Scientific Computing}, 39(1):A141--A163, 2017.

\bibitem{schwab2011sparse}
Christoph Schwab and Claude~Jeffrey Gittelson.
\newblock Sparse tensor discretizations of high-dimensional parametric and
  stochastic pdes.
\newblock {\em Acta Numerica}, 20:291--467, 2011.

\bibitem{ss11}
David Silvester and Valeria Simoncini.
\newblock An optimal iterative solver for symmetric indefinite systems stemming
  from mixed approximation.
\newblock {\em ACM Transactions on Mathematical Software}, 37(4), 2011.

\bibitem{Valeria}
Valeria Simoncini.
\newblock Computational methods for linear matrix equations.
\newblock {\em SIAM Review}, 58(3):377--441, 2016.

\bibitem{stoverud2016poro}
Karen~H St{\o}verud, Martin Aln{\ae}s, Hans~Petter Langtangen, Victor Haughton,
  and Kent-Andr{\'e} Mardal.
\newblock Poro-elastic modeling of syringomyelia--a systematic study of the
  effects of pia mater, central canal, median fissure, white and gray matter on
  pressure wave propagation and fluid movement within the cervical spinal cord.
\newblock {\em Computer methods in biomechanics and biomedical engineering},
  19(6):686--698, 2016.

\bibitem{wang2000theory}
Herbert~F Wang.
\newblock Theory of linear poroelasticity, 2000.

\end{thebibliography}

\end{document}